\documentclass[12pt]{article}
\usepackage{amsfonts,amsmath,amssymb,amscd,amsthm}
\usepackage[russian]{babel}
\def\Z{{\mathbb Z}} \def\R{{\mathbb R}}
\def\C{{\mathbb C}} \def\Q{{\mathbb Q}}

\long\def\comment#1\endcomment{}

\input{xy}
\xyoption{all}

\def\Z{{\mathbb Z}} \def\R{{\mathbb R}}
\def\C{{\mathbb C}} \def\Q{{\mathbb Q}}

\long\def\comment#1\endcomment{}

\newtheoremstyle{mydefinition}
  {3pt}
  {3pt}
  {\normalfont}
  {\parindent}
  {\bfseries}
  {.}
  { }
  {}

\theoremstyle{mydefinition}
\newtheorem{pr}{}[section]

\begin{document}



\comment


\smallskip
{\bf \ref{soproh}.}
Обозначим через $g$ данный многочлен.
По формуле косинуса тройного угла число $\cos(2\pi/9)$ является корнем многочлена $8x^3-6x+1$.
Значит, $\cos(2\pi/9)$ является корнем остатка $k+lx+mx^2$ от деления многочлена $g(1+x+x^2)$ на $8x^3-6x+1$ (здесь не предполагается, что $m\ne0$).
Если $k=l=m=n=0$, то $g(1+x+x^2)$ делится на $8x^3-6x+1$.
Так как  число $\cos(8\pi/9)$ также является корнем многочлена $8x^3-6x+1$, то
оно является корнем многочлена $g(1+x+x^2)$.

Приведем теперь к противоречию предположение о том, что многочлен $k+lx+mx^2$ ненулевой.
Если он ненулевой, то число $\cos(2\pi/9)$ является корнем остатка $ax+b$ от деления многочлена $8x^3-6x+1$ на $k+lx+mx^2$ (здесь не предполагается, что $a\ne0$).
По теореме о рациональных корнях уравнение $y^3-3y+1=0$ не имеет рациональных корней
(ибо числа $y=\pm 1$ не подходят).
Значит, уравнение $8x^3-6x+1=0$ не имеет рациональных корней.
Поэтому остаток $ax+b$ не может быть

 $\bullet$ нулевым, поскольку тогда $8x^3-6x+1$ делится на $k+lx+mx^2$, а значит, имеет множитель первой степени и, стало быть, рациональный корень,

 $\bullet$ ненулевым, ибо тогда число $\cos(2\pi/9)$ является корнем этого остатка, рационально, но это число является также корнем многочлена $8x^3-6x+1$.

Противоречие.

\endcomment

\centerline{\uppercase{\bf Some more proofs from the Book: }}
\bigskip
\centerline{\uppercase{\bf solvability and insolvability of equations in radicals}
\footnote{This note is based on lectures at summer school `Modern Mathematics' (Dubna), Kirov region summer school,  Moscow `olympic' school, Summer Conference of Tournament of Towns, as well as at circles `Mathematical seminar'
(Kolmogorov college) and `Olympiades and Mathematics' (MCCME).
I would like to acknowledge A.Ya. Belov, I.I. Bogdanov, G.R. Chelnokov, J. Cigler,
A.L. Glazman, A.S. Golovanov, A.A. Kaznacheev, P.V. Kozlov,
A.A. Pakharev, V.V. Prasolov, L.M. Samoilov, V.V. Volkov, M.N. Vyalyi and J. Zung for useful discussions.
}
}
\bigskip
\centerline{{\bf A. Skopenkov}
\footnote{skopenko@mccme.ru; Info: www.mccme.ru/\~{}skopenko;
Moscow Institute of Physics and Technology  and the Independent University of Moscow.
Supported by Simons-IUM fellowship. }
}

\small
\bigskip
{\bf Abstract.}
This paper is purely expository.
We present short elementary proofs of

$\bullet$ the Gauss Theorem on the constructibility of regular polygons by means of compass and ruler;

$\bullet$ the existence of a cubic equation unsolvable in real radicals;

$\bullet$ the Galois Theorem on the existence of a quintic equation unsolvable in complex radicals.

The statements of these celebrated results are simple and well-known.
However, their proofs given in most textbooks rely upon much unmotivated material and are far from being economic.
We do not use the terms `Galois group' or even `group'.
The paper is accessible for students familiar with
polynomials and complex numbers, and could be an interesting easy reading for professional mathematicians.

\bigskip
\normalsize
\hfill{\it Listeners are prepared to accept unstated (but hinted)}

\hfill{\it generalizations much more than they are able ... to}

\hfill{\it decode a precisely stated abstraction and to re-invent}

\hfill{\it the special cases that motivated it in the first place.}

\hfill{\it P. Halmos, How to talk mathematics.}

\subsection*{Statements of results}

A calculator (which calculates with absolute precision and has infinite memory) has operations
$$1\text{ (getting number 1)} ,\quad + \text{ (summation)},\quad - \text{ (subtraction)},$$
$$\times \text{ (multiplication)},
\quad : \text{ (division)}\quad\mbox{and}\quad \sqrt[n]{}
\text{ ($n$-power root)}.$$
A {\it real calculator} operates with real numbers and gives an error on attempting to extract
an even degree root of a negative number.
(As a more formal definition of a number that can be obtained on a real calculator a reader can accept the assertion of the Tower of Extensions Lemma below.)

A real number is {\it real constructible}, if it can be obtained on the real calculator so that all the roots extracted will be of degree 2.

\smallskip
{\bf The Gauss Theorem.}
{\it The number $\displaystyle\cos\frac{2\pi}n$ is real constructible if and only if $n=2^\alpha p_1\dots p_l,$ where $p_1,\dots,p_l$
are distinct primes of the form $2^{2^s}+1$.}

\smallskip
The following is apparently a folklore result.

\smallskip
{\bf Theorem on the insolvability in real radicals.}
{\it There exists a polynomial of the third degree with rational coefficients (e.g. $8x^3-6x+1$) none of whose roots can be obtained at a real calculator.}

A {\it complex calculator} has the same buttons as the real one, but operates with complex numbers and after pressing the button $\sqrt[n]{}$
gives all the $n$ values of the root.
A number is said to be obtainable on the complex calculator if a {\it set} of numbers containing given one can be obtained.
(This is not natural for a calculator; a reader can accept as a definition the assertion of the Tower of Extensions Lemma.)

\smallskip
{\bf The Galois Theorem.}
{\it There exists a polynomial of the fifth degree with rational coefficients (e.g. $x^5-4x+2$) none of whose roots can be obtained on the complex  calculator.}

\smallskip
This theorem is implied by the following result which is interesting and nontrivial  even for polynomials of degree 5.

\smallskip
{\bf The Kronecker Theorem.} {\it If a polynomial  with rational coefficients  is irreducible over $\Q$ and has
prime degree, more than one real root and at least one non-real root, then none of the roots can be obtained on
the complex calculator.}

\smallskip
In this note we present {\it short elementary proofs of these celebrated results}.
\footnote{Proof of the `if' part of the Gauss theorem is borrowed from [E1, \S24]; we simplify it by not using [E1, \S24, Lemma 2];
some motivation for the proof can be found in [E1].
Proof of the `only if' part of the Gauss theorem is based on [D, Supplement to \S\S35-37].
Proof of the Gauss theorem presented here
is simpler than the proof in
[KS];
proof of the `if' part from [KS] is added as an appendix.
Proof of the Galois Theorem is borrowed from [D, P, T]; we make certain simplifications and corrections.
This proof is different from the proof in [FT, Sk1].}
The term
`Galois group' (or even `group') is not used.
However, our presentation may be a good way to learn (or recall) starting ideas of the Galois theory, cf. [E1, E2, L, S, Sk].
For history of these results see [CR].

Let
$$\varepsilon_n:=\cos\dfrac{2\pi}n+i\sin\dfrac{2\pi}n.$$


\subsection*{A reduction of the Gauss theorem to complex numbers}

A complex number is {\it constructible}, if it can be obtained on the complex calculator so that all the roots extracted are of degree 2.


\smallskip
{\bf Complexification Lemma.} {\it A complex number is constructible if and only if its real and
imaginary parts are real-constructible.}

\smallskip
{\it Hint.} The `if' part is clear.
In order to prove the `only if' part write $\sqrt{a+bi}=u+vi$ and express $u,v$ by
quadratic radicals of $a$ and $b$. QED

\subsection*{Proof of the `if' part of the Gauss theorem }

 {\it Proof of the `if' part of the Gauss theorem for $n=5$.} (This proof is formally not required for the general case.)
Set
$$\varepsilon:=\varepsilon_5,\qquad T_0:=\varepsilon+\varepsilon^2+\varepsilon^4+\varepsilon^3=-1,$$
$$T_1:=\varepsilon+i\varepsilon^2-\varepsilon^4-i\varepsilon^3,
\qquad T_2:=\varepsilon-\varepsilon^2+\varepsilon^4-\varepsilon^3\qquad\text{and}
\qquad T_3:=\varepsilon-i\varepsilon^2-\varepsilon^4+i\varepsilon^3.$$
Then $\varepsilon=\dfrac{T_0+T_1+T_2+T_3}4$.
So by the Complexification Lemma it suffices to prove that $T_1,T_2,T_3$ are constructible.

Clearly, $T_1$ goes to $-iT_1$ under the substitution of $\varepsilon$ to $\varepsilon^2$.
Hence $T_1^4$ is invariant under this substitution.
Open brackets in the product $T_1^4$ and substitute $\varepsilon^5$ to 1.
We obtain a number represented in the form $a_0+a_1\varepsilon+a_2\varepsilon^2+a_3\varepsilon^3+a_4\varepsilon^4$ with $a_k\in\Z+i\Z$.
Since $T_1^4$ is invariant under substitution of $\varepsilon$ to $\varepsilon^2$, we have
$a_1=a_2=a_4=a_3$.
\footnote{This needs a justification, see the proof for the general case below or [E, Lemma 2 in p. 27].}
Therefore $T_1^4=a_0-a_1\in \Z+i\Z$.
Thus $T_1$ is constructible.
Analogously $T_2$ and $T_3$ are constructible.
QED

\smallskip
{\bf Multiplication Lemma.} {\it If $\varepsilon_m$ and $\varepsilon_n$
are constructible and $m,n$ are relatively prime, then $\varepsilon_{2n}$ and $\varepsilon_{mn}$ are constructible.}

\smallskip
{\it Hint.}
The constructibility of $\varepsilon_{2n}$ follows because $\varepsilon_{2n}=\pm\sqrt{\varepsilon_n}$.
The constructibility of $\varepsilon_{mn}$ follows because $\varepsilon_{mn}=\varepsilon_m^x\varepsilon_n^y$,
where $x$ and $y$ are integers such that $mx+ny=1$.
QED

\smallskip
{\bf Primitive Root Theorem.}
{\it For each prime $p$ there exists an integer $g$ (a} primitive root modulo {\it $p$) such that the residues modulo $p$ of $g^1,g^2,g^3\dots,g^{p-1}$ are distinct.}

\smallskip
{\it Hint for $p=2^m+1$ (only this case is used for the Gauss Theorem).}
If there are no primitive roots, then the congruence $x^{2^{m-1}}\equiv1\ (p)$ has $p-1=2^m>2^{m-1}$ solutions.
QED

\smallskip
{\it Proof of the `if' part of the Gauss theorem.}
If $l$ is odd, then $2^{kl}+1$ is divisible by $2^k+1$.
Thus if $2^m+1$ is a prime then $m$ is a power of 2.
So by the Complexification and the Multiplication Lemmas it suffices to prove that $\varepsilon_n$ is constructible for $n=2^m+1$ a prime.
Let $g$ be a primitive root modulo $n$.
 Since $n=2^m+1$, by the Multiplication Lemma $\beta:=\varepsilon_{n-1}$ is constructible.
For
$$r=0,1,2,\dots,n-2\quad\text{let}\quad T_r(x):=x+\beta^rx^g+\beta^{2r}x^{g^2}+\dots+\beta^{(n-2)r}x^{g^{n-2}}\in\Z[\beta][x].$$
(`Lagrange resolutions.')
Then $\varepsilon=\dfrac{(T_0+T_1+\dots+T_{n-2})(\varepsilon)}{n-1}$.
We have $T_0(\varepsilon)=-1$.
So it suffices to prove that $T_r(\varepsilon)$ is constructible for each $r=1,2,\dots,n-2$.

Since
$$T_r(x^g)\equiv\beta^{-r}T_r(x)\mod(x^n-1),\qquad\text{we have}
\qquad T_r^{n-1}(x^g)\equiv T_r^{n-1}(x)\mod(x^n-1).$$
Take the polynomial $a_0+a_1x+a_2x^2+\dots+a_{n-1}x^{n-1}$ with coefficients in $\Z[\beta]$, that is congruent to
$T_r^{n-1}(x)$ modulo $x^n-1$.
Then $a_k=a_{kg\mod n}$ for each $k=1,2,\dots,n-1$.
Hence $a_1=a_2=\dots=a_{n-1}$.
Therefore $T_r^{p-1}(\varepsilon)=a_0-a_1\in \Z[\beta]$.
So $T_r(\varepsilon)$ is constructible.
QED

\subsection*{Appendix: Gauss' effective proof of the `if' part}

This proof is longer than the above one, but gives a more realistic opportunity to obtain explicit formulas [BK, Sa].
(It would be also interesting to obtain explicit formulas using the above method.)

The material before the proof of the `if' part of the Gauss theorem is only required to show how the proof could be invented,
and could be omitted if a reader wants just a formal proof.

\smallskip
{\it Effective proof of the `if' part of the Gauss theorem for $n=5$.}
It suffices to prove that $\varepsilon:=\varepsilon_5$ is constructible.
We shall first prove that {\it some polynomials} of $\varepsilon$ are constructible.
Since
$$1+\varepsilon+\varepsilon^2+\varepsilon^3+\varepsilon^4=0,\quad\mbox{we have}
\quad(\varepsilon+\varepsilon^4)(\varepsilon^2+\varepsilon^3)=
\varepsilon+\varepsilon^2+\varepsilon^3+\varepsilon^4=-1.$$
Let
$$A_0:=\varepsilon+\varepsilon^4\quad\mbox{and}
\quad A_1:=\varepsilon^2+\varepsilon^3.$$
Then $A_0$ and $A_1$ are roots of the equation $t^2+t-1=0$ by the Vieta
theorem.
Hence these numbers are constructible.
Since $\varepsilon\cdot \varepsilon^4=1$, the numbers $\varepsilon$ and
$\varepsilon^4$ are roots of the equation
$t^2-A_0t+1=0$ by the Vieta theorem.
Therefore $\varepsilon$ (and $\varepsilon^4$) is constructible.

\smallskip
{\it Idea of Gauss' effective proof of the `if' part of the Gauss theorem.}
Assume that $n=2^m+1$.
It suffices to prove that $\varepsilon_n$ is constructible.
First it would be good to split the sum
$$\varepsilon+\varepsilon^2+\dots+\varepsilon^{n-1}=-1$$
into two summands $A_0$ and $A_1$ whose {\it product} is constructible
(or, in other words, to {\it group} the roots of the equation
$1+\varepsilon+\varepsilon^2+\dots+\varepsilon^{n-1}=0$ in a clever way).
Then $A_0$ and $A_1$ would be constructible by the Vieta Theorem.

Next it would be good to split the sum $A_0$ into two summands
$A_{00}$ and $A_{01}$ whose product is constructible, and analogously split $A_1=A_{10}+A_{11}$.
And so on, until we obtain that $A_{0\dots0}=\varepsilon$ is constructible.

It is however quite non-trivial to find the necessary splittings.
This is possible only for a {\it prime} $n=2^m+1$.
This is done using the Primitive Root Theorem.

\smallskip
{\it Gauss'  effective proof of the `if' part.}
If $l$ is odd, then $2^{kl}+1$ is divisible by $2^k+1$.
Thus if $2^m+1$ is a prime then $m$ is a power of 2.
So by the Complexification and the Multiplication Lemmas it suffices to prove that $\varepsilon_n$ is constructible for $n=2^m+1$ a prime.

Assume that $a_i\in\Z_2=\{0,1\}$ for each $i\in\{0,1,2,\dots,m-1\}$.
Let
$$\overline{a_{m-1}\dots a_1a_0}:=a_0+2a_1+2^2a_2+\dots+2^{m-1}.$$
(`Binary expansion.')
It is important for our proof that zeroes at the beginning are not thrown away.

Let $g$ be a primitive root modulo $n$.
For $A\in\Z_2^k$ let
$$T_A:=\sum_{B\in\Z_2^{m-k}} \varepsilon^{g^{\overline{BA}}}.$$
Let us prove by induction on $k$ that {\it for each $k$ and $A\in\Z_2^k$ the number $T_A$ is constructible.}
Then for $k=m$ we obtain that $T_{\overline{\underbrace{0\dots0}_m}}=\varepsilon$ is constructible.

The base $k=0$ follows because $T_\emptyset=-1$.
Let us prove the inductive step.
We have $T_A=T_{\overline{0A}}+T_{\overline{1A}}$.
Furthermore,
$$T_{0A}T_{1A}=\sum_{s=0}^{n-1}N_s\varepsilon_n^s\overset{(*)}=
N_0+\sum\limits_{C\in\Z_2^k}N_{g^{\overline C}}T_C.$$
Here the number $N_s$ (depending on $A$) is the quantity of ordered solutions
$$B_0,B_1\in\Z_2^{mk-1}\quad\text{of the congruence}\quad g^{\overline{B_00A}}+g^{\overline{B_11A}}\equiv s\mod n.$$
It is clear that $N_s=N_{g^{2^k}s}$.
Hence we have (*).
Hence $T_{0A}$ and $T_{1A}$ are expressed in terms of square radicals.
This implies the inductive step.
QED


\subsection*{Tower of extensions and the notion of a field}

If $F\subset\C$, $r\in\C$ и $r^q\in F$ for some positive integer $q$, then let
$$F[r]:=\{a_0+a_1r+a_2r^2+\dots+a_{q-1}r^{q-1}\ |\ a_0,\dots,a_{q-1}\in F\}.$$

{\bf Tower of Extensions Lemma.} {\it Let $F\in\{\R,\C\}$.
A number $x\in F$ can be obtained at the $F$-calculator if and only if there are
$r_1,\ldots r_{s-1}\in F$ and $q_1,\ldots q_{s-1}\in\Z$ such that $q_k\ge 2$ and
$$\Q=F_1\subset F_2\subset F_3\subset \ldots\subset F_{s-1}\subset F_s\ni x,
\quad\mbox{where}\quad r_k^{q_k}\in F_k,\quad r_k\not\in F_k\quad\mbox{and}\quad F_{k+1}=F_k[r_k]$$
for each $k=1,\dots,s-1$.}

\smallskip
{\it Hint.} Use induction on the number of operations of the calculator,
which are necessary to construct given number $x$.
In the
inductive step use multiplication by conjugate.
QED

\smallskip
Such a sequence is called {\it a tower of extensions}.

A {\it field} is a subset of $\C$ which is closed under summation, subtraction, multiplication and division by a non-zero number.
This notion is useful for us because theorem on division with a remainder and its corollaries hold for polynomials with coefficients in a field.

If $F$ is a field, $q$ a prime and  $r\not\in F$, then the polynomial $t^q-r^q$ is irreducible over $F$
(this is trivial for $q=2$ and essentially proved in the Linear Independence Lemma below for arbitrary $q$).
Hence $F[r]$ is a field.

\subsection*{Proof of the `only if' part of the Gauss theorem }

Since $\varepsilon_n=\varepsilon_{nk}^k$, the constructibility of $\varepsilon_{nk}$ implies the constructibility of $\varepsilon_n$.
Hence by the Complexification Lemma in order to prove the `only if' part of the Gauss theorem we need
to prove that $\varepsilon_n$ is not constructible for

(A) $n$ a prime not of the form $2^m+1$, and

(B) $n=p^2$ the square of a prime.


\smallskip
{\bf Lemma on Powers of 2.}
{\it If a polynomial with coefficients in $\Q$ is irreducible over $\Q$ and has a constructible root, then the degree of the polynomial
is a power of 2.}

\smallskip
This lemma is proved below.

Now the non-constructibility follows by applying the Tower of Extensions Lemma and the Lemma on Powers of 2 to
the root $\varepsilon_n$ of the polynomial

$\bullet$ $P(x):=x^{n-1}+x^{n-2}+\dots+x+1$ for case (A) and

$\bullet$ $P(x):=x^{p(p-1)}+x^{p(p-2)}+\dots+x^p+1$ for case (B).

The irreducibility of $P(x)$ over $\Z$ follows (in both cases) by the irreducibility of $P(x+1)$ over $\Z$.
The latter is implied by the following Eisenstein criterion.

{\it Let $p$ be a prime. If the
leading coefficient of a polynomial with integer coefficients is
not divisible by $p$, other coefficients are divisible by $p$ and
the constant term is not divisible by $p^2$, then this polynomial
is irreducible over $\Z$.}

The irreducibility over $\Q$ follows by the irreducibility over $\Z$ and the following
Gauss lemma.

{\it If a polynomial with integer
coefficients is irreducible over $\Z$, then it is irreducible over $\Q$.}

Both the Eisenstein criterion and the Gauss lemma are easily proved by passing to polynomials with coefficients in
$\Z_p$ (for the Gauss lemma take a decomposition $P=P_1P_2$ of given polynomial $P$ over $\Q$,
take $N_1$ and $N_2$ such that both
$N_1P_1$ and $N_2P_2$ have integer  coefficients and take a prime divisor $p$ of $N_1N_2$).

\smallskip
{\bf Conjugation Lemma.} {\it Let $F\subset\C$ be a field, $r\in\C-F$ and $r^2\in F$.
Define the conjugation map
$$\overline\cdot:Q_{k+1}[\sqrt{a_k}]\to Q_{k+1}[\sqrt{a_k}]\quad\mbox{by}
\quad\overline{x+y\sqrt{a_k}}:=x-y\sqrt{a_k}.$$
This map is well-defined,}
$$\overline{z+w}=\overline z+\overline w, \quad \overline{zw}=\overline z\cdot\overline w\quad\mbox{and} \quad
\overline z=z\Leftrightarrow z=x+0\sqrt{a_k}\in Q_k.$$

The Lemma on Powers of 2 is the case $k=1$ of the following assertion.

\smallskip
{\bf Generalized Lemma on Powers of 2.} {\it If
$$\Q=F_1\subset F_2\subset F_3\subset \ldots\subset F_{s-1}\subset F_s\ni \alpha,
\quad\mbox{where}\quad r_k^2\in F_k,\quad r_k\not\in F_k\quad\mbox{and}\quad F_{k+1}=F_k[r_k]$$
for each $k=1,\dots,s-1$.
If $P$ is a polynomial with coefficients in $F_k$ such that $P(\alpha)=0$ and $P$ is irreducible over $F_k$, then $\deg P$ is a power of 2.  }

\smallskip
{\it Proof.}
Downward induction on $k$.
The base $k=s$ is clear.
Let us prove the inductive step.
Denote by $P_k$ any polynomial with coefficients in $F_k$ irreducible over $F_k$ and having the root $\alpha$.
Consider divisibility and GCD in $F_{k+1}$.
Since $P_k(\alpha)=0$, the polynomial $P_k$ is divisible by $P_{k+1}$.
By the Conjugation Lemma for $F=F_k$ and $F[r]=F_{k+1}$ we have $P_k=\overline{P_k}$ is divisible by $\overline{P_{k+1}}$.
Let $D:=GCD(P_{k+1},\overline{P_{k+1}})$.
Since $P_{k+1}$ is irreducible over $F_{k+1}$ and divisible by $D$, it follows that either $D=1$ or $P_{k+1}=D$.

In the second case since $\overline D=D$, we have $P_{k+1}=D\in F_k[x]$.
Hence $P_k=P_{k+1}$ and the inductive step is proved.

In the first case $P_k$ is divisible by $M:=P_{k+1}\overline{P_{k+1}}$.
Since $\overline M =M$, we have $M\in F_k[x]$.
Since $P_k$ is irreducible over $F_k$, we have $P_k=M$.
Hence $\deg P_k=2\deg P_{k+1}$ is a power of 2 by the inductive hypothesis.
QED

\subsection*{Proof of the insolvability in real radicals}

{\bf Main Lemma (real case).} {\it Suppose that $q$ is a prime, $F\subset \R$ a field, $r\in\R-F$ and $r^q\in F$.

(a) (linear independence)
If $P(r)=0$ for some polynomial $P\in F[\varepsilon_q][t]$ of degree less than $q$, then $P=0$.

Or, equivalently, for each element $\alpha\in F[\varepsilon_q,r]$ there are unique
$a_0,a_1,\dots,a_{q-1}\in F[\varepsilon_q]$ such that $x=a_0+a_1r+a_2r^2+\dots+a_{q-1}r^{q-1}$.

(b) (conjugation) If $P\in F[\varepsilon_q][t]$ and $P(r)=0$, then $P(r\varepsilon_q^k)=0$ for each $k=0,1,\dots,q-1$.}

\smallskip
{\it Proof of (a).}
Both polynomials $P$ and $t^q-r^q$ with coefficients in $F[\varepsilon_q]$ have the root $r$.
Hence their GCD has a root $r$ and degree $k$, $0<k\le\deg P<q$.
All the roots of the polynomial $t^q-r^q$ are $r,r\varepsilon_q,r\varepsilon_q^2, \dots, r\varepsilon_q^{q-1}$.
The free coefficient of the GCD is the product of some $k$ of these roots.
Then $r^k\in F[\varepsilon_q]$.
Since $q$ is prime, $kx+qy=1 $ for some integers $x,y$.
Then $r=(r^k)^x(r^q)^y\in F[\varepsilon_q]$.
\footnote{Using the fact that dimension is well-defined the next paragraph can be rewritten as $\dim_FF[r]\le\dim_F F[\varepsilon_q]\le q-1$.}

Therefore $r^2,r^3,\dots,r^{q-1}\in F[\varepsilon_q]$.
Take a table of size $q\times(q-1)$ with entries $a_{kl}\in F$ formed by representations of numbers $r^k$ in powers of $\varepsilon_q$:
$$r^k = \sum \limits_{l=0}^{q-2}a_{kl}\varepsilon_q^l, \qquad 0 \le k \le q-1. $$
Using some operations of addition of one line multiplied by a number from $F$ to another line,
we can obtain a table with a zero line.

Hence there is a nonzero polynomial $P_1$ of degree less than $q$ with coefficients in $F$ and the root $r$.
Further arguments are similar to the first paragraph.
It is only necessary to replace $P$ by $P_1$ and $F[\varepsilon_q]$ by $F$.
We obtain that $r\in F$.
A contradiction.
QED

\smallskip
{\it Proof of (b).}
Since $P(r)=0$, the remainder of division of $P(t)$ by $t^q-r^q$ assumes the value 0 at $r$.
Hence by (a) this remainder is zero.
This implies the conclusion of (b).
QED

\smallskip
{\it Proof of the insolvability in real radicals.}
Suppose to the contrary that a root $x_0$ of the equation  $8x^3-6x+1=0$ can be obtained on the real calculator.
Then by the Tower of Extensions Lemma for $F=\R$ there is the smallest $s$ for which there is a tower of extensions whose
last field $F_s$ contains a root $x_1$ of the equation  $8x^3-6x+1=0$ (possibly, $x_1\ne x_0$).
Let $F:=F_{s-1}$, $q:=q_{s-1}$ and $r:=r_ {s-1}$.
Then $x_1=h(r)$ for some polynomial $h$ with coefficients in $F$ of degree greater than 0 and less than $q$.

Apply the Main Lemma (real case) (b) to $P(t):=8h(t)^3-6h(t)+1$.
Since $8h(r)^3-6h(r)+1=0$, we obtain that $h(r\varepsilon_q^k)$ is the root of $8x^3-6x+1=0$ for each $k=0,1,\dots,q-1$.
If $h(r\varepsilon_q^k)=h(r\varepsilon_q^l)$ for some $0\le k<l\le q-1$, then by the Main Lemma (real case) (a)
$\deg h=0$.
This contradicts the minimality of $s$.

Thus the numbers $h(r\varepsilon_q^k)$, $0\le k\le q-1$, are distinct roots of the equation $8x^3-6x+1=0$.
So $q=2$ or $q=3$.
If $q=2$, then by the Vieta theorem the third root of the equation $8x^3-6x+1=0$ equals $-2h(0)\in F$.
Contradiction with the minimality of $s$.
If $q=3$, then $\overline{\varepsilon_3} = \varepsilon_3^2$ implies $\overline{h(r\varepsilon_3)}=h(r\varepsilon_3^2)$.
This contradicts the fact that the last two numbers are real and distinct.
QED

\subsection*{Proof of the Kronecker Theorem}

Of the following lemmas only the Main Lemma (a,c) (linear independence and decomposition) and the
Lemma on an Improvement of Tower of Extensions are directly used in the proof of the Kronecker Theorem.
Main Lemma (b) (conjugation) is only used for Main Lemma (c) (decomposition);
the Gauss Lowering Lemma is only used for the Lemma on an Improvement of Tower of Extensions.

Main Lemma (c) (decomposition) shows that under certain conditions all the roots of a polynomial $g\in F[x]$
are values of a polynomial from $F[x]$ at conjugate points.

\smallskip
{\bf Main Lemma.} {\it Suppose that $q$ is a prime, $F\subset \C$ a field, $r\in\C-F$ and $r^q, \varepsilon_q\in F$.

(a) (linear independence)
If $P(r)=0$ for a polynomial $P\in F[t]$ of degree less than $q$, then $P=0$.

Or, equivalently, for each $\alpha\in F[\varepsilon_q,r]$ there are unique
$a_0,a_1,\dots,a_{q-1}\in F[\varepsilon_q]$ such that $x=a_0+a_1r+a_2r^2+\dots+a_{q-1}r^{q-1}$.

(b) (conjugation) If $P\in F[x,t]$ and $P(x,r)=0$ as a polynomial in $x$, then $P(x,r\varepsilon_q^k)=0$ as a polynomial in $x$ for each $k=0,1,\dots,q-1$.

(c) (decomposition) If the polynomial $g$ of prime degree with coefficients in $F$ is irreducible over $F$ and
is reducible over $F[r]$, then $g(x)=A(x-x_0)(x-x_1)\dots(x-x_ {q-1})$ for some $A\in F$ and distinct values $x_0,x_1,\dots,x_{q-1}$ of a polynomial with coefficients in $F$ at points $r,r\varepsilon_q,\dots,r\varepsilon_q^{q-1}$.}

\smallskip
{\it Proof of (a).}
Part (a) is analogous to the first paragraph of the proof of the Main Lemma (real case) (a) in \S7, only replace $F[\varepsilon_q]$ by $F$.
QED

\smallskip
{\it Proof of (b).}
The statement of (b) is invariant under the division of the polynomial $P$ by $t^q-r^q$ with a remainder.
Therefore we can assume that $\deg_tP<q$.
In this case (b) is obtained by applying (a) to the coefficients.
QED

\smallskip
{\it Proof of (c).}
By the assumption there is a irreducible (over $F[r]$) divisor of the polynomial $g$ in $F[r]$.
We may assume that the top degree coefficient of this divisor is 1.
This divisor is the value of $h(x,r)$ at $r$ of a polynomial $h\in F[x,t]$ of degrees in more than 0 in $t$ and less than $\deg g$ in $x$.
So $h(x,r)$ is irreducible over $F[r]$ and $g(x)=h(x,r)h_1(x,r)$ for some polynomial $h_1\in F[x,t]$.
Let $\varepsilon:=\varepsilon_q$.

We apply (b) to $P(x,r):=g(x)-h(x,r)h_1(x,r)$.
We obtain that $g(x)$ is divisible by the polynomial $h(x,r\varepsilon^k)$ in $F[r]$ for each $k=0,1,\dots,q-1$.

The polynomial $h(x,r\varepsilon^k) $ is irreducible over $F[r]$ for each $k=0,1,\dots,q-1$.

(Otherwise apply (b) to the polynomial $P$ that is the difference of $h(x,r\varepsilon^k)$ and its factors.
We obtain that the polynomial $h(x,r)$ is reducible over $F[r]$.
Contradiction.)

By (a) the polynomials $h(x,r\varepsilon^k)$ are different for different $k=0,1,\dots,q-1$.
Hence $g$ is divisible by their product.
By (a) this product can be uniquely represented in the form
$$a_0(x) + a_1(x)r + \dots + a_{q-1}(x)r^{q-1}\quad\text{for some}\quad a_k\in F[x].$$
Since the product goes to itself under the change $r\to r\varepsilon$ (which is well defined by (a)),
by (a) $a_k(x)=a_k(x)\varepsilon^k\in F[x]$ for each $k=1,2,\dots,q-1$.
Hence $a_k(x)=0$ for each $k=1,2,\dots,q-1$.
This means that the product is $a_0(x)\in F[x]$.

From this and the irreducibility of $g$ over $F$ it follows that $g$ equals to this product.
Then $\deg g=q\deg_xh$.
Since $\deg g$ is a prime and $\deg_xh<\deg g$, we have $\deg_xh=1$ (and $\deg g=q$).
So, $-h(0,r)\in F[r]$ is a root of $g$.
And the remaining roots are $-h(0,r\varepsilon^k)$ for $k=0,1,\dots,q-1$.
QED

\smallskip
{\bf Gauss Lowering Lemma.}
{\it For each $n$ the number $\varepsilon_n$ can be obtained on the complex calculator so that all the roots extracted have degrees strictly less than $n$.}

\smallskip
{\it Proof.} By induction on $n$.
The base is obvious.
If $n=ab$ for some integers $0<a,b<n$, then the induction step follows from $\varepsilon_n=\sqrt[a]\varepsilon_b$.
If $n$ is a prime, then the inductive step follows by the following assertion whose proof is analogous to the
proof of constructibility in the Gauss theorem.

{\it If $n$ is prime, then the $\varepsilon_n$ can be obtained on the complex calculator so that all the roots extracted have degree $n-1$.}
QED

\smallskip
{\bf Lemma on an Improvement of Tower of Extensions.}
{\it If a number can be obtained on the complex calculator, then there is a tower of extensions such that for each
$k=1,2,\dots,s-1$ the number $q_k$ is prime, $\varepsilon_{q_k}\in F_k$, and either $r_k\in\R$ or $|r_k|^2\in F_k$.}

\smallskip
{\it Proof.}
By downward induction on $q$ let us show that from an arbitrary tower of extensions one can obtain a tower of extensions for which $\varepsilon_{q_k}\in F_k$ for each $q_k>q$.
If we show that, for $q=1$ we obtain a tower of extensions for which $\varepsilon_{q_k}\in F_k$ for each $k=1,2,\dots,s-1$.

The base: $q=\max_kq_k$; in this case there is nothing to prove.
To prove the inductive step we take the smallest $k$ for which $q_k=q$.
Paste between $F_{k-1}$ and $F_k$ `getting $\varepsilon_q$
extracting only roots of degree less than $q$' obtained from the Gauss Lowering Lemma.
If such $k$ does  not exist, then the inductive step is obvious.

Next, replace each extraction of a root of a composite degree $ab$ by extraction of roots of $a$-th and $b$-th degrees.
The condition $\varepsilon_{q_k}\in F_k$ for each $k=1,2,\dots,s-1$ is preserved, because if $\varepsilon_{ab}\in F_k$,
then $\varepsilon_a\in F_k$ and $\varepsilon_b\in F_k$.

Call a tower of extensions {\it interesting} if for each $k=1,2,\dots,s-1$ the number $q_k$ is prime and $\varepsilon_{q_k}\in F_k$.
By downward induction on $l$ let us show that from an arbitrary interesting tower of extensions one can obtain an interesting tower of extensions
such that for each $k\le s-l$
$$\overline F_k=F_k\quad\text{and \quad either}\quad r_k\in\R \quad\text{or}\quad|r_k|^2\in F_k.$$
Then for $l=0$ we obtain the Lemma.
Base: $l=s-1$, in which case there is nothing to prove.
Let us prove the inductive step.
(If $r_k\in\R$, then the inductive step is obvious, but the following argument also works.)
Since $F_k=\overline{F_k}$ and $r_k^{q_k}\in F_k$, we obtain $|r_k|^{2q_k}=r_k^{q_k}\overline{r_k^{q_k}}\in F_k$.
So $F_k[|r_k|^2]=F_k[\sqrt[q_k]{|r_k|^{2q_k}}]$, where the real value of the root is taken.
Replace subtower
$$F_k\subset F_{k +1}\subset\dots\subset F_s\quad\text{by subtower}\quad F_k \subset F_k[|r_k|^2]
\subset F_k[r_k,\overline r_k]=F_{k+1}[\overline r_k]\subset\dots\subset F_s[\overline r_k].$$
Clearly, the property of being interesting is preserved under this change.
In the new subtower replace repeated copies of the same field by one field.
After that we apply the inductive hypothesis.
QED

\smallskip
{\it Proof of the Kronecker Theorem.}
Suppose to the contrary that a root of the polynomial $g$ can be obtained on the complex calculator.
Then take a tower of extensions from the Lemma on an Improvement of Tower of Extensions.
Since $g$ is irreducible over $\Q$ and reducible over the last field of the tower, there is $s$ such that $g$ is irreducible over $F_s$ and reducible over $F_{s+1}$.
Let $r:=r_s$ and $q:=q_s$.
By the Main Lemma (c) $g(x)=A(x-x_0)(x-x_1)\dots(x-x_{q-1})$ for some $A\in F_s$ and different values $x_0,x_1,\dots,x_{q-1}$
of a polynomial $a_0+a_1t+\dots+a_{q-1}t^{q-1}$ with coefficients in $F_s$ at points $r\varepsilon_q^k$, $0\le k\le q-1$.
The property $x_k\in\R$ is equivalent to $x_k=\overline x_k$.
Note that $\overline{\varepsilon_q^k}=\varepsilon_q^{-k}$.

If $r\in\R$, then by the Main Lemma (a) for each $k\in\{0,1,\dots,q-1\}$ the condition $x_k=\overline x_k$ is equivalent to
$a_s\varepsilon^{2sk}=\overline a_s$ for each $s=0,1,\dots,q-1$.
Consequently, $x_k\in\R$ for at most one $k$.

If $r\not\in\R$, then by the Improvement Lemma $|r|^2\in F_s$.
Then $\overline r^s=\dfrac{|r|^{2s}}{r^q}r^{qs}$, where $\dfrac{|r|^{2s}}{r^q}\in F_s$.
Hence, by the Main Lemma (a) for each $k\in\{0,1,\dots,q-1\}$ the condition $x_k=\overline x_k$ is equivalent
to `$a_0=\overline a_0$ and $a_s=\overline a_{qs}\dfrac{|r|^{2q-2s}}{r^q}$ for each $s=1,2,\dots,q-1$'.
These equations do not depend on $k$.
Therefore if one of the numbers $x_0,\dots,x_{q-1}$ is real, then they are all real.

Contradiction.
QED

\bigskip
\newpage
\centerline{\bf References}

[B] J. Bergen, A Concrete Approach to Abstract Algebra: From the Integers to the
Insolvability of the Quintic, 2010.
http://www.elsevierdirect.com/product.jsp?isbn=9780123749413

[BK] Burda Yu. and Kadets L., Seventeen-sides-polygon and the Gauss reciprocity law (in Russian), Mat. Prosveschenie, 17 (2013).
http://www.mccme.ru/free-books/matprosi.html.

[CR] R. Courant and H. Robbins, What is Mathematics, Oxford Univ. Press.

[D] H. D\"orrie, 100 Great Problems of Elementary Mathematics: Their History and Solution, Dover Publ., New York, 1965.

[E1] H.M. Edwards, Galois Theory, Springer Verlag, 1984.

[E2] H.M. Edwards, The construction of solvable polynomials,
Bull. Amer. Math. Soc. 46 (2009), 397-411.
Errata: Bull. Amer. Math. Soc. 46 (2009), 703-704.

http://www.ams.org/journals/bull/2009-46-03/S0273-0979-09-01253-1/

[FT] D. Fuchs, S. Tabachnikov,  Mathematical Omnibus. AMS, 2007.
\linebreak
http://www.math.psu.edu/tabachni/Books/taba.pdf


[H] Ch.R. Hadlock, Field Theory and its Classical Problems, Carus Mathematical Monographs 19,
The Mathematical Association of America, 1978.
\linebreak
http://books.google.com/books?id=5s1p0CyafnEC\&printsec=frontcover\&dq=hadlock.

[KS] P. Kozlov and A. Skopenkov, A la recherche de l'alg\`ebre perdue:
du cote de chez Gauss, Mat. Prosveschenie, 12 (2008), 127--144. arxiv:math/0804.4357

[Le] L. Lerner, Galois Theory without abstract algebra, http://arxiv.org/abs/1108.4593.

[P] V.V. Prasolov, Problems in algebra, arithmetics and analysis, Moscpw, MCCME, 2007,
ftp://ftp.mccme.ru/users/prasolov/algebra/algebra2.pdf

[S] J. Stillwell, Galois theory for beginners, Amer. Math. Monthly, 101 (1994), 22-27.

[Sk] A. Skopenkov, Philosophical and methodical appendix, in: Mathematics as  a sequence of problems, Ed. A. Zaslavsky, D. Permyakov, A. Skopenkov, M. Skopenkov, A. Shapovalov, MCCME, Moscow, 2009.

[Sk1] A. Skopenkov, A simple proof of the Abel-Ruffini theorem,
Mat. Prosveschenie, 15 (2011) 113-126,
http://arxiv.org/abs/1102.2100.

[Sa] A. Safin, A program for construction of regular polygons by compass and ruler,
\linebreak
http://www.mccme.ru/mmks/dec08/Safin.pdf

[T] V.M. Tikhomirov, Abel and his great theorem, Kvant, 2003, N1.
\linebreak
http://kvant.mccme.ru/pdf/2003/01/kv0103abel.pdf

\newpage
\normalsize
\centerline{\uppercase{\bf Еще несколько доказательств из Книги: }}
\bigskip
\centerline{\uppercase{\bf разрешимость и неразрешимость уравнений в радикалах}
\footnote{
Заметка основана на занятиях, проведенных в ЛШ `Современная математика', Кировской ЛМШ, Московской ВШ,
а также на кружках `Математический семинар' и `Олимпиады и математика'.
Благодарю А.Я. Белова-Канеля, И.И. Богданова, Э.Б. Винберга, В.В. Волкова, М.Н. Вялого, А.С. Голованова, П.А. Дергача, Д. Зунга, А.А. Казначеева,
В.А. Клепцына, Г.А. Мерзона, А.А. Пахарева, В.В. Прасолова, А.Д. Руховича, Л.М. Самойлова, М.Б. Скопенкова, Г.Р. Челнокова, Л.А. Шабанова и
В.В. Шувалова за полезные
замечания и предложения.}
}
\bigskip
\centerline{\bf А. Скопенков
\footnote{Поддержан грантом фонда Саймонса.
Московский Физико-Технический Институт,
Независимый Московский Университет,
Инфо: www.mccme.ru/\~{ }skopenko}
}

\bigskip
\small
{\bf Аннотация.}
Основное содержание этой брошюры --- простые  элементарные доказательства знаменитых теорем Гаусса,
Абеля, Галуа и Кронекера о построимости правильных многоугольников и неразрешимости уравнений в радикалах.
На примере этих доказательств иллюстрируются некоторые основные идеи алгебры.
Определения построимости и разрешимости в радикалах приводится;
для понимания доказательств достаточно знакомства с многочленами и умения извлекать корни из комплексных чисел.
Брошюра адресована всем любителям изложения глубоких идей на примерах красивых результатов и доказательств:
старшекласникам,
студентам,
учителям,
и профессиональным математикам.

\normalsize

\tableofcontents

\section{Введение}

\subsection{О чем эта брошюра}

Основное содержание этой брошюры --- простые элементарные доказательства

$\bullet$ теоремы Гаусса о построимости правильных многоугольников
(и даже более сильного результата --- теоремы Гаусса о понижении, \S\ref{gaucon});

$\bullet$ существования уравнения 3-й степени, неразрешимого в {\it вещественных} радикалах
(и даже более сильного результата --- сильной вещественной теоремы о неразрешимости, \S\ref{prostr});

$\bullet$ теоремы Галуа о существовании уравнения 5-й степени, неразрешимого в {\it комплексных} радикалах
(и даже более сильных результатов --- теоремы Кронекера, \S\ref{progal}).

Определения построимости и разрешимости в радикалах, а также формулировки указанных теорем, приведены в \S\S\ref{intsqu},\ref{intrad}.
Где найти их доказательства и вообще как устроена брошюра, написано в \S\ref{intplan}.

Брошюра адресована тем, кому интересен хотя бы один из этих результатов;
она может быть занимательным чтением для профессиональных математиков (им достаточно прочитать \S\ref{pro}).

Приводимые доказательства не претендуют на новизну (конкретные ссылки и приведены в в \S\S\ref{intsqu},\ref{intrad}).
Однако, к сожалению, приводимые доказательства малоизвестны.
Как следствие, малоизвестно, что не только решать квадратные и кубические уравнения, но и доказывать указанные
теоремы экономнее напрямую, а не строя и затем применяя теорию Галуа, ср. [Kh2, Ki].
При этом, конечно, на таких прямых доказательствах ясно видны базовые идеи теории Галуа
(см. подробнее \S\ref{motphi}).

Доказательство разрешимости основано на методе {\it резольвент Лагранжа};
доказательства неразрешимости основаны на идее {\it сопряжения}, ср. [Va].
(В доказательстве сильной вещественной теоремы о неразрешимости из \S\ref{prostr} используется также теорема о размерности башни расширений.)
Разбор доказательств (или их начала) полезен для закрепления тем `многочлены', `комплексные числа',
`иррациональность' и `основы линейной алгебры' на кружке или в матклассе.
Старшеклассники найдут в заметке задачи для иследования, не претендующие на научную новизну, см. подробнее \S\ref{motint}.

\subsection{Разрешимость в квадратных радикалах: формулировки}\label{intsqu}

 Известно, что
$$\cos\frac{2\pi}3=-\frac12,\quad \cos\frac{2\pi}4=0, \quad \cos\frac{2\pi}5=\frac{\sqrt5-1}4, \quad \cos\frac{2\pi}6=\frac12, \quad \cos\frac{2\pi}8=\frac1{\sqrt2}.$$
Как обобщить эти формулы (используя только четыре арифметические действия и извлечения корней)?

Рассмотрим калькулятор с кнопками
$$1,\quad +,\quad -,\quad \times,\quad :\quad\mbox{и}\quad \sqrt[n]{}\quad
\mbox{для любого $n$}.$$
Калькулятор вычисляет числа с абсолютной точностью и имеет неограниченную память.
При делении на 0 он выдает ошибку.

Пусть сначала калькулятор {\it вещественный}, т.е.  оперирует с вещественными числами
и при извлечении корня четной степени из отрицательного числа выдает ошибку.

Вещественное число называется {\it вещественно построимым}, если его можно получить на вещественном калькуляторе
так, чтобы при этом извлекались корни только второй степени (т.е. получить из 1 при помощи сложений, вычитаний, умножений, делений и извлечений квадратного корня из положительных чисел).

Например, вещественно построимы числа
$$
\sqrt[4]2=\sqrt{\sqrt2},\quad \sqrt{2\sqrt3},\quad \sqrt2+\sqrt3,\quad
\sqrt{1+\sqrt2},\quad1+\sqrt{3-2\sqrt2},\quad
\frac1{1+\sqrt2}\quad\mbox{и}\quad \cos 3^\circ.$$
Про последнее число это не совсем очевидно (но доказано в \S\ref{gauref}).

\smallskip
{\bf Теорема Гаусса.} {\it Число $\cos(2\pi/n)$ вещественно построимо тогда и только тогда,
когда $n=2^\alpha p_1\dots p_l$, где $p_1,\dots,p_l$ --- различные простые числа вида $2^{2^s}+1$.}

\smallskip
Теорему Гаусса можно переформулировать в терминах {\it построимости циркулем и линейкой} правильных многоугольников.
См. \S\ref{motcon}; эта переформулировка не используется в остальном тексте.

Строго говоря, теорема Гаусса не дает настоящего решения проблемы построимости, поскольку неизвестно, какие числа
вида $2^{2^s}+1$ являются простыми.
Однако теорема Гаусса дает, например, быстрый
алгоритм выяснения построимости числа $\cos(2\pi/n)$.

История этой знаменитой теоремы здесь не приводится, см. [Gi].
\footnote{Доказательство построимости в теореме Гаусса получено из [E1, \S24]
некоторым упрощением (мы обходим использование леммы 2).
Оно более простое по сравнению с доказательством из [KS].
Элементарное доказательство построимости для $n=17$ приводится, например, в [BK, Ch, D, \S37, Gi, P, Po, PS, Ko]
(при этом иногда приводятся явные формулы, как с доказательствами утверждений о знаках перед радикалами
[D, \S37, Sa], так и без [BK]).
Для общего случая оно намечено в [Ga, Gi], где ясности доказательства немного
мешает построение общей теории вместо доказательства конкретного результата.
Подход из [K] дает объяснение на вопрос `почему', и было бы интересно довести его до полного доказательства.
\newline
Доказательство непостроимости в теореме Гаусса основано на [D, Supplement to \S\S35-37].
Оно более простое по сравнению с доказательствами из [KS].
\newline
Ср. [H, B, Vi, W].
См. подробнее \S\ref{motint}. }

Вещественная непостроимость числа $\cos(2\pi/9)$ влечет следующий результат.

\smallskip
{\bf Следствие.}  {\it Трисекция угла невозможна на вещественном калькуляторе, если можно извлекать корни только второй степени.
Или, формально, число $\cos(\alpha/3)$ невозможно получить на нем, имея число $\cos\alpha$ (например, для $\alpha=2\pi/3$).}

\smallskip
{\bf Замечание.}
Для практики приближенные методы вычисления тригонометрических функций и решения уравнений
более полезны, чем радикальные формулы.
Кроме того, уравнения степени выше 4 разрешимы при помощи трансцедентных функций (см. метод Виета в \S\ref{mot34} и [PS]).
Однако проблема разрешимости в радикалах интересна как пробная задача современной теории {\it исследовании операций}.

\subsection{Неразрешимость в радикалах: формулировки}\label{intrad}

Следующее утверждение хорошо известно, см. \S\ref{mot34}.

{\it Если многочлен третьей степени с рациональными коэффициентами имеет ровно один вещественный корень,
то этот корень можно получить на вещественном калькуляторе.
\footnote{Стандартная терминология: уравнение {\it разрешимо в вещественных радикалах}.}
 Более того, это можно сделать так, чтобы извлечение корня происходило
только два раза, один раз второй и один раз третьей степени. }

Следующая теорема, по-видимому, является фольклорным результатом.

\smallskip
{\bf Теорема о неразрешимости в вещественных радикалах.}
{\it Существует многочлен 3-й степени с рациональными коэффициентами (например, $8x^3-6x+1$),
ни один из корней которого невозможно получить на вещественном калькуляторе.}

\smallskip
{\bf Следствие.} {\it Трисекция угла невозможна на вещественном калькуляторе.
Или, формально, число $\cos(\alpha/3)$ невозможно получить на нем, имея число $\cos\alpha$
(например, для $\alpha=2\pi/3$). }

\smallskip
{\it Доказательство.}
По формуле косинуса тройного угла каждое из чисел $\cos(2\pi/9)$, $\cos(8\pi/9)$, $\cos(14\pi/9)$ удовлетворяет
уравнению  $8x^3-6x+1=0$.
Значит, то теореме ни одно из них невозможно получить на вещественном калькуляторе.
QED

\smallskip

{\it Комплексный} калькулятор имеет те же кнопки, что и вещественный, но оперирует с комплексными числами и при нажатии
кнопки $\sqrt[n]{}$ выдает все значения корня.
На комплексном калькуляторе можно получить число, если на нем можно получить {\it множество} чисел, содержащих
заданное число.

Следующее утверждение хорошо известно, см. \S\ref{mot34}.

{\it Все корни любого многочлена третьей или четвертой степени с рациональными коэффициентами
можно получить на комплексном калькуляторе.
\footnote{Стандартная терминология: уравнение {\it разрешимо в радикалах}.}
Более того, это можно сделать так, чтобы извлечение корня происходило только

$\bullet$ два раза, причем один раз третьей степени и один раз второй --- для многочлена третьей степени.

$\bullet$ четыре раза, причем один раз третьей степени и три раза второй --- для многочлена четвертой степени.}

\smallskip
{\bf Теорема Галуа.}
{\it Существует многочлен 5-й степени (например, $x^5-4x+2$), ни один из корней которого невозможно получить на комплексном калькуляторе.
\footnote{Немного ранее была доказана более слабая теорема П. Руффини - Н.Х. Абеля.
Она сложнее формулируется [A, FT, S], но более знаменита, ибо именно она решила знаменитую проблему о
разрешимости уравнений в радикалах.
Наиболее простое известное мне доказательство теоремы Руффини-Абеля (приводимое здесь) дает сразу теорему Галуа.}
}

\smallskip
Из приведенных теорем о неразрешимости легко вывести, что для любого $n\ge3$ ($n\ge5$) существует многочлен $n$-й степени, ни один из корней которого невозможно получить на вещественном (комплексном) калькуляторе.
 \footnote{Доказательство теорем о неразрешимости в вещественных и комплексных радикалах
основано на замечательных статье [T] и книгах  [D, \S25, P, дополнение 8]
(впрочем, здесь исправлены неточности, см. \S\ref{progal}).
Ср. [Sa1] (этот текст написан Л. Самойловым по курсу, проведенному совместно В. Волковым и автором).
Оно отлично от доказательства из [A, FT, S].
Почему корни любого уравнения степени ниже 5 выражаются в радикалах через коэффициенты, а степени 5 и выше --- нет?
Доказательство  из [A, FT, S] дает такой ответ: поскольку группа $S_n$ разрешима в точности при $n\le4$.
Приводимое доказательство дает такой ответ: поскольку 5 простое и больше 3.
Простота `причины' неразрешимости косвенно указывает на простоту доказательства.
См. подробнее \S\ref{motint}. }

\subsection{План брошюры}\label{intplan}

Эту брошюру не обязательно читать подряд.
Читатель может выбрать удобную ему последовательность изучения (или вовсе опустить некоторые пункты) на основании приводимого плана.
К нему разумно вернуться, если читатель потеряет нить изложения.

В \S\S\ref{motint},\ref{motphi} приводятся мотивировки.
В \S\S\ref{motcon} приводится переформулировка теоремы Гаусса.
В \S\ref{mot34} приводится решение уравнений 3-й и 4-й степени (в задачах);
вместо немотивированных замен мы естественно сводим уравнения к таким, которые ясно, как решать.
Пункты в \S\ref{mot} независимы с остальным текстом
(т.е. ни один из них не используются в остальном тексте и для изучения любого из них достаточно прочитать \S1).

Вот схема зависимости остальных пунктов.
$$\xymatrix{& & \ref{nonn2} \ar@{-->}[rrd] \ar[r] &  \ref{nongau} \ar@{-->}[rrd] &  \\
& \ref{non12} \ar[ur] \ar[r] \ar[dr] & \ref{non14} \ar@{-->}[u] \ar@{-->}[d] &  & \ref{profie} \ar[r] \ar[rd]\ar[rdd] & \ref{progau} \ar@{-->}[d]\\
\ref{gauref} \ar[dr] &\ref{gaueff} & \ref{non13} \ar[r] & \ref{non1p} \ar[r] & \ref{nonn} \ar@{-->}[u] \ar@{-->}[r] & \ref{prorea} \ar@{-->}[d]\\
\ref{gaures} \ar[r] & \ref{gaucon} \ar@{-->}[u] \ar[rrr] & & & \ref{progal} \ar[r] & \ref{prostr}
}$$
Пунктир в схеме означает, что один параграф нужен для мотивировки другого, но формально не используется в другом.
Если в одном пункте используется малая явно указанная часть другого, то это не считается формальной зависимостью.

В \S\S\ref{gauref},\ref{gaucon} доказана построимость в теореме Гаусса.
В \S\ref{gaures} основные идеи этого доказательства иллюстрируются в задачах.
В \S\ref{gaueff} приводится дополнительный материал.


В \S\ref{non} приведены задачи, подводящие к доказательствам неразрешимости
(т.е. непостроимости в теореме Гаусса теоремам о неразрешимости в радикалах).
В \S\ref{hints} приведены их решения.
При этом к доказательству непостроимости в теореме Гаусса (\S\ref{progau}) подводят только задачи из \S\S\ref{non12}--\ref{nongau},
а к доказательству теорем о неразрешимости в радикалах (\S\S\ref{prorea},\ref{progal}) --- задачи из \S\S\ref{non12}, \ref{non13}--\ref{nonn}.

Сами доказательства неразрешимости приводятся в \S\ref{pro}.
Пункты \S\ref{progau}--\ref{progal} формально независимы друг от друга.
Неформально же, для каждого следующего пункта полезно прочитать предыдущий.
В \S\ref{prostr} приводится дополнительный материал.


\section{Мотивировки и отступления}\label{mot}

\subsection{Чем интересны приводимые доказательства}\label{motint}

Приводимые доказательства намного проще и короче тех, которые излагаются в стандартных учебниках по алгебре.
(Здесь я имею в виду доказательства `с нуля', а не вывод нужной теоремы из построенной перед этим теории, в которой фактически заключается все доказательство.
Ср. с последней сноской в \S\ref{intrad} по поводу менее стандартных изложений.)
Так получилось благодаря тому, что в отличие от большинства учебников ,
приводимые доказательства не используют термина `группа Галуа' (даже термина `группа').
Несмотря на отсутствие этих {\it терминов}, {\it идеи} приводимых доказательств являются {\it отправными}
для теории Галуа и {\it конструктивной теории Галуа} [E2].
Ср. [Ch1, Le].
Более подробно это обсуждается в философско-методическом отступлении (\S\ref{motphi}).

Основные идеи представлены на `олимпиадных' примерах: на простейших частных случаях,
свободных от технических деталей, и со сведением к необходимому минимуму алгебраического языка.

Перед доказательством непостроимости в теореме Гаусса (\S\ref{progau}) некоторые его идеи
демонстрируются по одной и на простейших примерах (\S\S\ref{non12}--\ref{nongau}).
Эти примеры доказывают неразрешимость классических задач древности об удвоении куба и трисекции угла.

Неразрешимость доказывается сначала при условии, что {\it корень извлекался только один раз} (в задачах в \S\S\ref{non12},\ref{non13},\ref{non1p}).
Благодаря этому основные идеи преподносятся на примере рациональных чисел
(а не произвольных полей и даже не полей из башни расширений).
Идея сопряжения и другие идеи заключены в леммах, выделенных в \S\S\ref{non},\S\ref{pro}.

Мы показываем, {\it как можно придумать} приводимые доказательства.
Пути к ним намечены в виде задач (\S\ref{gaures} и \S\ref{non}).
Это характерно не только для дзенских монастырей, но и для серьезного преподавания математики.
К задачам приведены указания и решения.
Хотя {\it придумать} доказательства непросто, {\it изложить} их можно коротко, см. \S\ref{gaucon} и \S\ref{pro}.
Освобождение доказательства от деталей, возникших при его придумывании и не нужных для него самого ---
важная часть его проверки.

Многие из приведенных задач --- удачные темы для исследовательских работ школьников, ср. [Ak].
Например,  \ref{sopr}.c, \ref{foex}.c, \ref{utv4}, \ref{twoexa}.h, \ref{3ex}.d, \ref{utv3}, \ref{kvadiffreal}.d, \ref{kronint1} (попроще),
\ref{tower2}.d, \ref{isslreal}.a (посложнее),
\ref{quairr}.bc, \ref{que3}.4, \ref{isslreal}.b, \ref{issl4} (решение мне неизвестно, но может оказаться несложным).
Большинство этих задач не претендуют на научную новизну.

{\it Общие замечания к формулировкам задач.}
Задачи обозначаются жирными цифрами.
Если условие задачи является утверждением, то в задаче требуется это утверждение доказать.
Если некоторая задача не получается, то читайте дальше --- соседние задачи могут оказаться подсказками.
(На занятии задача-подсказка выдается только тогда, когда школьник  или студент подумал над самой задачей.)
Некоторые задачи не нужны для подведения к доказательству сформулированных теорем, но интересны сами по себе.

\subsection{Философско-методическое отступление}\label{motphi}


По моему мнению, именно с {\it новых идей\/}, изложенных на уже имеющемся языке, а не с
{\it введения нового языка\/}, полезно {\it начинать} изучение любой теории.
Как правило, такие идеи наиболее ярко выражаются доказательствами, подобными приведенным здесь.

{\it При изложении материала нужно ориентироваться на объекты, которые
основательнее всего укореняются в человеческой памяти.
Это --- отнюдь не системы аксиом и не логические приемы в доказательстве
теорем.
Изящное решение красивой задачи, формулировка которой ясна и доступна, имеет
больше шансов удержаться в памяти студента, нежели абстрактная теория.
Скажем больше, именно по такому решению, при наличии некоторой математической
культуры, студент впоследствии сможет восстановить теоретический материал.
Обратное же, как показывает опыт, практически невозможно} [Ko, предисловие].

Известно также, что `путь познания должен повторять путь развития'.
\footnote{Впрочем, это не вполне верно. Так, изучение геометрии
Лобачевского вовсе не обязательно начинать с попыток доказать Пятый
Постулат. Геометрия Лобачевского для нас сейчас важна, в первую
очередь, ее приложениями в ТФКП, теории чисел, топологии, теории
групп, алгебраической геометрии, космологии и т.д., а вовсе не тем,
что она демонстрирует независимость Пятого Постулата от остальных
аксиом Евклида. С этой точки зрения более плодотворно ее построение
не на основе аксиом Евклида-Гильберта, а на основе понятия группы
преобразований (Клейн) или римановой метрики (Риман). Аналогично,
изучение теории Галуа вовсе не обязательно начинать с задачи о решении
алгебраического уравнения в радикалах или квадратных радикалах.
С современной точки зрения теория Галуа есть теория алгебраических
расширений полей, составляющая неотъемлемую часть алгебры и имеющая
приложения и аналоги в других разделах математики (алгебраическая
геометрия, теория накрытий, теория инвариантов), а решение
алгебраических уравнений в радикалах --- это маргинальная задача.
(Э. Б. Винберг).}

Такой стиль изложения не только делает материал более доступным,
но позволяет сильным студентам (для которых доступно даже абстрактное
изложение) приобрести математический вкус и стиль с тем, чтобы

(1) разумно выбирать проблемы для исследования и их мотивировки.
\footnote{Математик, понимающий, что теория Галуа мотивируется более важными проблемами,
чем построимость правильных многоугольников и разрешимость
алгебраических уравнений в радикалах, вряд ли станет мотивировать созданную им
теорию приложениями, которые можно получить и без его теории.}

(2) ясно излагать собственные открытия, не скрывая ошибки или известности
полученного результата за чрезмерным формализмом.
\footnote{К сожалению, такое --- обычно бессознательное --- сокрытие ошибки часто
происходит с молодыми математиками, воспитанными на чрезмерно формальных курсах.
Происходило и с автором этих строк; к счастью, все мои серьезные ошибки исправлялись {\it перед} публикациями.}

Мода на искусственно формализованное изложение
\footnote{Видимо, общепринятый термин `бурбакизация' не очень удачен ввиду
`масштаба и влияния деятельности Бурбаки, независимо от оценки пользы и вреда разных ее аспектов' (А. Шень).}
привела к следующему парадоксу.
По данному {\it известному понятию} высшей математики зачастую непросто
восстановить {\it конкретный красивый результат}, для которого это понятие действительно необходимо
(и при получении которого это понятие возникло).

Доказательство c использованием некоторого нового термина имеют свои
преимущества:
оно подготавливает читателя к доказательству тех теорем, которые уже трудно или
невозможно доказать без этого термина.
\footnote{Например, векторное доказательство теоремы Пифагора уже является
достаточным основанием для введения понятий векторного пространства и
скалярного умножения, хотя эти понятия и не являются необходимыми для
доказательства упомянутой теоремы. (Э. Б. Винберг.) }
Однако такие доказательства, как правило, не должны быть {\it первыми}
доказательствами данного результата (легко себе представить результат
{\it первого} знакомства с теоремой Пифагора на основе понятий векторного
пространства и скалярного умножения).
Кроме того, при приведении `терминологического' доказательства полезно
оговорить его мотивированность не доказываемым результатом, а обучением полезному новому методу (ср. с (1) выше).

Приведенная выше точка зрения разделяется многими математиками (а некоторыми --- нет); я унаследовал ее от Ю. П. Соловьева.

Приводимые порой в качестве {\it основных} приложений теории Галуа
теоремы Гаусса и о неразрешимости уравнений в радикалах
неубедительны для мотивировки этой теории (так же, как приложение к решению
квадратных уравнений неубедительно для мотивировки общей теории разрешимости
уравнений произвольной степени в радикалах).
 Действительно, эти теоремы имеют элементарное доказательство, не использующее `групп Галуа'.
В терминах теории Галуа формулируется общий критерий разрешимости
алгебраического уравнения в радикалах.
Но этот критерий не дает настоящего решения проблемы разрешимости, а лишь сводит ее к трудной
задаче вычисления группы Галуа уравнения.
(То, что никакая {\it другая теория} не дает легкого для применений ответа,
не позволяет утверждать, что {\it теория Галуа} дает такой ответ.)
Но, конечно, формулировка общего критерия в адекватных проблеме терминах
может иметь важное философское значение.

Однако теория Галуа выходит далеко за рамки проблемы разрешимости уравнений в радикалах.
Ее популяризации послужит дальнейшая публикация интересных теорем, формулируемых без понятий
теории Галуа, но при попытках доказать которые она естественно возникает.
Примеры таких теорем мне сообщили А.Я. Белов, С.М. Львовский и Г.Р. Челноков
(к сожалению, в доступной мне начальной учебной литературе по теории Галуа мне
не удалось найти такие теоремы, формулировка которых не была бы скрыта под
толщей обозначений и терминов).


\subsection{Отступление: связь с построениями циркулем и линейкой}\label{motcon}

Используя отрезки длины $a$, $b$ и $c$, можно построить
циркулем и линейкой отрезки длины $a+b$, $a-b$, $ab/c$, $\sqrt{ab}$.
Поэтому если на плоскости задан отрезок длины 1, то отрезок вещественно построимой длины
можно построить циркулем и линейкой.
Этот простой результат (известный еще древним грекам) показывает, что из вещественной построимости числа $\cos(2\pi/n)$
вытекает построимость  правильного $n$-угольника циркулем и линейкой.

\smallskip
{\bf Основная теорема теории геометрических построений.}
{\it
Если отрезок длины $a$ можно построить
циркулем и линейкой, имея отрезок длины 1, то число $a$ вещественно построимо.}

\smallskip
Этот несложный результат [P, Ko] (доказанный лишь в 19-м веке)
показывает, что из непостроимости числа $\cos(2\pi/n)$ вытекает
непостроимость правильного $n$-угольника циркулем и линейкой.
Для доказательства этого результата можно рассмотреть все возможные случаи появления новых
объектов (точек, прямых, окружностей) и показать, что координаты всех построенных точек и
коэффициенты уравнений всех проведенных прямых и окружностей являются построимыми.
См. детали в [Ko, CR, Ma, P].


\subsection{Отступление: решение уравнений 3-й и 4-й степени}\label{mot34}

\begin{pr}\label{2zamena}
 (a) Уравнение $ax^3+bx^2+cx+d=0$ сводится к~уравнению $x^3+px+q=0$ заменой переменной.

(b) Уравнение $ax^4+bx^3+cx^2+dx+e=0$ сводится к~уравнению $x^4+px^2+qx+r=0$ заменой переменной.
\end{pr}

\begin{pr}\label{2center}
 (a) Найдите координаты центра симметрии графика функции $y=-2x^3-6x^2+4$.

(b) График любого кубического многочлена имеет центр симметрии.
\end{pr}

\begin{pr}\label{2numroo}
  Сколько корней имеет уравнение \qquad

(a) $x^3+2x+7=0$? \qquad (b) $x^3-3x+1=0$?
\end{pr}

В этой и следующей задачах можно пользоваться без доказательства теоремой о промежуточных значениях многочлена.

\begin{pr}\label{2numbroot}
  Найдите количество решений уравнения $x^3+px+q=0$ (в зависимости от параметров $p,q$).
\end{pr}

 \begin{pr}\label{2=1}
$\sqrt[3]{2+\sqrt5}-\sqrt[3]{\sqrt5-2}=1$.
\end{pr}

Некоторые из вышеприведенных задач сложны для начинающего.
Следующие задачи (из второй серии) являются подсказками.

\begin{pr}\label{2cubic}
 (a) Найдите хотя бы одно решение уравнения $x^3-3\sqrt[3]2x+3=0$. \quad

{\it Указание}. Так как $(b+c)^3=b^3+c^3+3bc(b+c)$, то число $b+c$ является корнем уравнения $x^3-3bcx-(b^3+c^3)=0$.

(b) Решите уравнение $x^3-3\sqrt[3]2x+3=0$.
\end{pr}


\begin{pr}\label{2ineq}
 (a) $a^2+b^2+c^2\ge ab+bc+ca$. Когда достигается равенство?
\qquad

(b) $a^3+b^3+c^3\ge3abc$ при $a,b,c>0$.
\end{pr}


\begin{pr}\label{2decom}
 (a) Разложите на множители выражение $a^3+b^3+c^3-3abc$.

(b) Разложите выражение $a^3+b^3+c^3-3abc$ на линейные множители с комплексными коэффициентами.
\end{pr}

\begin{pr}\label{2dferro}
  {\it Метод дель Ферро.}
Напишите формулу для решения уравнения $x^3+px+q=0$ предложенным методом.
При каком условии применим этот метод, если квадратные корни
разрешается извлекать только из положительных чисел?
\end{pr}

\begin{pr}\label{2calc}
  Калькулятор имеет кнопки $1,\quad +,\quad -,\quad \times\quad\mbox{и}\quad :\quad$.
Калькулятор вычисляет числа с абсолютной точностью и имеет неограниченную память.
При делении на 0
он выдает ошибку.

(a) Докажите, что любое кубическое уравнение можно решить на калькуляторе, если
дополнительно разрешается
два раза извлечь корень из положительного числа.

(b)*  Какие кубические уравнения можно решить на калькуляторе, если
дополнительно разрешается
один раз извлечь корень из положительного числа?
\end{pr}

\begin{pr}\label{2fourth}
  Решите уравнение \quad

(a) $(x^2+2)^2=18(x-1)^2$. \quad
(b) $x^4+4x-1=0$. \quad
(c) $x^3-3x-1=0$.
\end{pr}

{\it Указание к \ref{2fourth}.b.} Подберите такие $\alpha,b,c$, чтобы
$x^4+4x-1=(x^2+\alpha)^2-(bx+c)^2$. Для этого найдите хотя бы
одно $\alpha$, для которого квадратный трехчлен
$(x^2+\alpha)^2-(x^4+4x-1)$ является полным квадратом.



\begin{pr}\label{2ferrari}
  {\it Метод Феррари.}
(a) Решите уравнение $x^4+2x^2-8x-4=0$.

(b) Найдите формулу для корней уравнения $x^4+px^2+qx+r=0$ методом, намеченным в задачах \ref{2fourth}.ab,
использующую корень $\alpha$ вспомогательного кубического уравнения.
Не забудьте разобрать все случаи!
\end{pr}


\begin{pr}\label{2vieta}
  {\it Метод Виета.} (a)
$\cos(\alpha+\beta)=\cos\alpha\cos\beta-\sin\alpha\sin\beta$ и
$\sin(\alpha+\beta)=\sin\alpha\cos\beta+\cos\alpha\sin\beta$.

(b) $\sin3\alpha=3\sin\alpha-4\sin^3\alpha$ и
$\cos3\alpha=4\cos^3\alpha-3\cos\alpha$.

(c) Решите уравнение $4x^3-3x=\frac12$.


(d) Используя функции $\cos$ и~$\arccos$, напишите общую формулу для
решения уравнения $x^3+px+q=0$ методом, намеченным в этой задаче.
При каком условии уравнение $x^3+px+q=0$ решается этим методом?
\end{pr}

\subsection{Указания и решения к задачам}

{\bf \ref{2zamena}.} Воспользуйтесь заменой переменной $y:=x+\frac b{3a}$ или $y:=x+\frac b{4a}$.

\smallskip
{\bf \ref{2numroo}.} (a) {\it Ответ}: 1. Так как $f(-2)<0$ и $f(1)>0$, то по  теореме о
промежуточных значениях многочлена корень имеется.
Ввиду монотонности корень только один.

\smallskip
{\bf \ref{2cubic}.}  {\it Ответ}: $x=-1-\sqrt[3]{2}$.

(a) \textit{Указание.} $x^3-3\sqrt[3]{2}x+3=x^3-3bcx+(b^3+c^3)$, где $b=1$, $c=\sqrt[3]{2}$.

(b) В силу задачи \ref{2decom}.a уравнение $x^3-3\sqrt[3]2x+3=0$ равносильно
уравнению
$$(x+b+c)(x^2+b^2+c^2-bc-bx-cx)=0\quad\text{с}\quad~b=1\quad\text{и}\quad и~c=\sqrt[3]2.$$
По задаче \ref{2ineq}.a второй сомножитель положителен при любом~$x$ (поскольку $b\ne c$).
Значит, исходное уравнение имеет единственное решение $x=-b-c=-1-\sqrt[3]2$.

\smallskip
{\bf \ref{2decom}.} (a) Поделите $a^3-3abc+(b^3+c^3)$ на $a+b+c$ `уголком'.

(b) {\it Ответ:} для $\varepsilon=(-1+i\sqrt 3)/2$
$$a^3+b^3+c^3-3abc=(a+b+c)(a^2+b^2+c^2-ab-bc-ca)=(a+b+c)(a+b\varepsilon+c\varepsilon^2)
(a+b\varepsilon^2+c\varepsilon).$$

{\bf \ref{2fourth}.} {\it Ответы}:  (a) $(-3\sqrt2\pm\sqrt{10+12\sqrt2})/2$;
\quad (b) $(-\sqrt2\pm\sqrt{4\sqrt2-2})/2$.

(с) {\it Указание}. Сведите к~\ref{2vieta}.c заменой $y=2x$.

{\it Ответ}: $2\cos\frac\pi9$, $2\cos\frac{7\pi}9$,
$2\cos\frac{13\pi}9=2\cos\frac{5\pi}9$.

\smallskip
{\it Замечание к \ref{2ferrari}.b.} Уравнение $x^4+ax^3+bx^2+cx+d=0$ можно также решить,
подобрав такие $\alpha$, $A$, $B$, что
$x^4+ax^3+bx^2+cx+d=\Big(x^2+\frac{ax}2+\alpha\Big)^2-(Ax+B)^2.$

\smallskip
{\bf \ref{2vieta}.} (c) {\it Ответ}: $\cos\frac\pi9$, $\cos\frac{7\pi}9$,
$\cos\frac{13\pi}9=\cos\frac{5\pi}9$.

\section{Доказательство построимости в теореме Гаусса}\label{gau}

\subsection{Переформулировка построимости в теореме Гаусса}\label{gauref}

\begin{pr}
Число $\cos(2\pi/n)$ вещественно построимо для $n=3,4,5,6,8,10,15$.
\end{pr}

\begin{pr}
{\bf Лемма об умножении (вещественная версия).}

(a) Если $\cos(2\pi/n)$ вещественно построимо, то $\cos(\pi/n)$ вещественно построимо.

(b) Если $\cos(2\pi/n)$ и $\cos(2\pi/m)$ вещественно построимы и $m,n$ взаимно просты, то $\cos(2\pi/mn)$ вещественно построимо.
\end{pr}

Из этой леммы вытекает, что вещественная построимость в теореме Гаусса следует из вещественной построимости
чисел $\cos(2\pi/n)$ для простых $n$ вида $2^{2^s}+1$.

Комплексное число называется {\it построимым}, если его можно получить на комплексном калькуляторе так,
чтобы при этом извлекались корни только второй степени.

\begin{pr}\label{g-eps}
Число $\cos(2\pi/n)$ построимо тогда и только тогда,
когда число
$$\varepsilon_n:=\cos(2\pi/n)+i\sin(2\pi/n)$$
построимо.
\end{pr}

\begin{pr}\label{g-com}
{\bf Лемма о комплексификации.}
Комплексное число построимо тогда и только тогда, когда
его вещественная и мнимая части вещественно построимы.
\end{pr}

Из этой леммы вытекает, что вещественное число построимо тогда и только тогда, когда оно вещественно построимо.
\footnote{Заметим, что на комплексном калькуляторе нет кнопок $Re$ и $Im$.
Однако их можно `реализовать', доказав, что если можно получить число $z$, то можно получить и $\overline z$.
Но так будет доказана {\it построимость} вещественной и мнимой части, а не их {\it вещественная построимость.}
Для доказательства вещественной построимости нужно научиться извлекать корень из комплексного числа при помощи вещественного калькулятора.
Это возможно только для корней второй степени.
Если в определении построимости и вещественной построимости допускать извлечения корней третьей степени, то
аналог леммы о комплексификации будет неверен (ибо существует уравнение третьей степени, разрешимое на комплексном калькуляторе, но не разрешимое на вещественном).}
Значит, построимость в теореме Гаусса достаточно доказать с заменой `вещественной построимости' на `построимость'.


\subsection{Метод резольвент Лагранжа}\label{gaures}

\begin{pr}  \label{postr7}
(a) Число $\varepsilon_5$ построимо.

(b) На комплексном калькуляторе можно получить число $\varepsilon_7$ так, чтобы при этом корни извлекались
только второй и третьей степени.

(b') На комплексном калькуляторе можно получить число $\varepsilon_7$ так, чтобы при этом корни второй и третьей степени извлекались только по одному разу, а корни большей степени не извлекались совсем.

(c) На комплексном калькуляторе можно получить число $\varepsilon_{11}$ так, чтобы при этом корни извлекались
только второй и пятой степени.

(d) Число $\varepsilon_{17}$ построимо.

(e) Докажите построимость в теореме Гаусса.
\end{pr}

Пункты (a), (b) и (b') можно просто решить непосредственно.
Для решения пунктов (c,d,e) уже нужны более продвинутые идеи, к изложению которых мы переходим.

\begin{pr}\label{ident} Решите системы уравнений ($x,y,z,t$ --- неизвестные, $a,b,c,d$ --- известные):
$$(a) \begin{cases}x+y+z+t=a\\
x+iy-z-it=b\\
x-y+z-t=c\\
x-iy-z+it=d
\end{cases}.
\qquad
(b) \begin{cases}x+y+z=a\\
x+\varepsilon y+\varepsilon^2z=b\\
x+\varepsilon^2y+\varepsilon z=c
\end{cases}, \quad\text{где}\quad \varepsilon:=\dfrac{-1+i\sqrt3}2.$$
\end{pr}

{\it Доказательство построимости числа $\varepsilon:=\varepsilon_5$.}
Обозначим
 \begin{align*}
T_0:=&\varepsilon+\varepsilon^2+\varepsilon^4+\varepsilon^8=-1,\\
T_1:=&\varepsilon+i\varepsilon^2-\varepsilon^4-i\varepsilon^8,\\
T_2:=&\varepsilon-\varepsilon^2+\varepsilon^4-\varepsilon^8\qquad\text{и}\\
T_3:=&\varepsilon-i\varepsilon^2-\varepsilon^4+i\varepsilon^8.
\end{align*}
(`Резольвенты Лагранжа'.)
 Тогда  $T_0+T_1+T_2+T_3=4\varepsilon$.
Поэтому достаточно доказать построимость каждого из чисел $T_1,T_2,T_3$.

При замене $\varepsilon$ на $\varepsilon^2$ число $T_1$ переходит в $-iT_1$.
Значит, $T_1^4$ не меняется при этой замене.
Поэтому $T_1^4$ не меняется при двукратной и трехкратной таких заменах, т.е., при заменах
$\varepsilon$ на $\varepsilon^4$ и $\varepsilon$ на $\varepsilon^8=\varepsilon^3$.
Итак, для любого $k$ число $T_1^4$ не меняется при замене $\varepsilon$ на $\varepsilon^k$.

Раскроем скобки в произведении $T_1^4$ и заменим $\varepsilon^5$ на 1.
Получим равенство
$$T_1^4=a_0+a_1\varepsilon+a_2\varepsilon^2+a_3\varepsilon^3+a_4\varepsilon^4\quad\text{для некоторых}\quad a_k\in\Z+i\Z.$$
Так как для любого $k$ число $T_1^4$ не меняется при замене $\varepsilon$ на $\varepsilon^k$, получаем $a_1=a_2=a_3=a_4$.
Поэтому $T_1^4=a_0-a_1\in \Z+i\Z$.
Значит, $T_1$ построимо.
Аналогично $T_2$ и $T_3$ построимы.

В приведенном рассуждении нужно обосновать вывод `$a_1=a_2=a_3=a_4$'
(и строго определить, что такое `замена $\varepsilon$ на $\varepsilon^2$').
Обоснование для общего случая трудное [E1, \S24].
Вместо этого обоснования изменим немного доказательство.
(Именно этим изменением приводимое доказательство отличается от данного в [E1, \S24].)

Определим многочлен $T_1(x):=x+ix^2-x^4-ix^8$.
Определим многочлены $T_0(x)$, $T_2(x)$ и $T_3(x)$ формулами, аналогичными вышенаписанным.
Как и выше, $(T_0+T_1+T_2+T_3)(\varepsilon)=4\varepsilon$.
Поэтому достаточно доказать построимость каждого из чисел $T_r(\varepsilon)$, $r=1,2,3$.
Имеем
\footnote{Два многочлена называются {\it сравнимыми по модулю многочлена $x^5-1$}, если их разность делится на $x^5-1$.}
$$iT_1(x^2)\equiv T_1(x)\mod(x^5-1)\quad\Rightarrow\quad T_1^4(x^2)\equiv T_1^4(x)\mod(x^5-1)\quad\Rightarrow$$
$$\Rightarrow\quad T_1^4(x^k)\equiv T_1^4(x)\mod(x^5-1)\quad\text{для любого }k.$$
Возьмем многочлен $a_0+a_1x+a_2x^2+a_3x^3+a_4x^4$ с коэффициентами в $\Z+i\Z$,
сравнимый с  $T_1^4(x)$ по модулю $x^5-1$.

 Тогда $a_1=a_2=a_3=a_4$.
 Поэтому $T_1^4(\varepsilon)=a_0-a_1\in \Z+i\Z$.
\footnote{Другой способ, предложенный М. Ягудиным:
$$T_1^4(\varepsilon)=a_0+a_1\varepsilon+a_2\varepsilon^2+a_3\varepsilon^3+a_4\varepsilon^4=
a_0+a_1\varepsilon^2+a_2\varepsilon^4+a_3\varepsilon+a_4\varepsilon^3=$$
$$=a_0+a_1\varepsilon^3+a_2\varepsilon+a_3\varepsilon^4+a_4\varepsilon^2=
a_0+a_1\varepsilon^4+a_2\varepsilon^3+a_3\varepsilon^2+a_4\varepsilon.$$
Суммируя эти выражения, получим $4T_1^4(\varepsilon)=a_0-a_1-a_2-a_3-a_4\in \Z+i\Z$.}

Значит, $T_1(\varepsilon)$ построимо.
Аналогично $T_2(\varepsilon)$ и $T_3(\varepsilon)$ построимы.
QED

\begin{pr}\label{7}
(a) Обозначим
$$\beta:=\dfrac{1+i\sqrt3}2\quad\text{и}
\quad T(x):=x+\beta x^3+\beta^2x^9+\beta^3x^{27}+\beta^4x^{81}+\beta^5x^{243}.$$
\quad
Докажите, что $T(x)\equiv\beta T(x^3)\mod(x^7-1)$.

(b) Обозначим
$$\beta:=\varepsilon_{10}\quad\text{и}
\quad T(x):=x+\beta x^2+\beta^2x^4+\beta^3x^8+\beta^4x^{16}+\dots+\beta^9x^{512}.$$
\quad
Докажите, что $T(x)\equiv\beta T(x^2)\mod(x^{11}-1)$.
\end{pr}

\subsection{Доказательство построимости в теореме Гаусса}\label{gaucon}

{\bf Лемма об умножении.} {\it (a) Если $\varepsilon_n$ построимо, то $\varepsilon_{2n}$ построимо.

(b) Если $\varepsilon_n$ и $\varepsilon_m$ построимы и $m,n$ взаимно просты, то $\varepsilon_{mn}$ построимо.}

\smallskip
{\it Доказательство} получается из формул $\varepsilon_{2n}\in\sqrt{\varepsilon_n}$ и
$\varepsilon_{mn}=\varepsilon_m^x\varepsilon_n^y$, где $x$ и $y$ --- целые числа, для которых $nx+my=1$.
QED

\smallskip
{\bf Теорема о первообразном корне.} {\it Для любого простого $p$ существует число $g$,
для которого остатки от деления на $p$ чисел $g^1,g^2,g^3,\dots,g^{p-1}$ различны.}

 \smallskip
{\it Указание к доказательству для $p=2^m+1$ (только этот случай нужен для теоремы Гаусса).}
Если первообразного корня нет, то сравнение $x^{2^{m-1}}\equiv1\bmod p$ имеет $p-1=2^m>2^{m-1}$ решений.
QED

\smallskip
{\it Доказательство построимости в теореме Гаусса.}
По леммам о комплексификации и об умножении достаточно доказать, что $\varepsilon_n$ построимо для простого $n=2^{2^s}+1$.
Так как $n-1=2^m$, то по лемме об умножении $\beta:=\varepsilon_{n-1}$ построимо.
Обозначим
$$\Z[\beta]:=\{a_0+a_1\beta+a_2\beta^2+\dots+a_{n-2}\beta^{n-2}\ |\ a_0,\dots,a_{n-2}\in\Z\}.$$
Пусть $g$ --- первообразный корень по модулю $n$.
 Для
$$r=0,1,2,\dots,n-2,\quad \text{обозначим}\quad
T_r(x):=x+\beta^rx^g+\beta^{2r}x^{g^2}+\dots+\beta^{(n-2)r}x^{g^{n-2}}\in\Z[\beta][x].$$
 Тогда
$(T_0+T_1+\dots+T_{n-2})(\varepsilon)=(n-1)\varepsilon.$
 Имеем $T_0(\varepsilon)=-1$.
Поэтому достаточно доказать построимость каждого из чисел $T_r(\varepsilon)$, $r=1,2,\dots,n-2$.
Имеем
$$\beta^rT_r(x^g)\equiv T_r(x)\mod(x^n-1)\quad\Rightarrow\quad T_r^{n-1}(x^g)\equiv T_r^{n-1}(x)\mod(x^n-1)\quad\Rightarrow$$
$$\Rightarrow\quad T_r^{n-1}(x^k)\equiv T_r^{n-1}(x)\mod(x^n-1)\quad\text{для любого }k.$$
Возьмем многочлен $a_0+a_1x+a_2x^2+\dots+a_{n-1}x^{n-1}$
с коэффициентами в $\Z[\beta]$,  сравнимый с  $T_r^{n-1}(x)$ по модулю $x^n-1$.
Тогда
$a_1=a_2=\dots=a_{n-1}$.
 Поэтому $T_r^{n-1}(\varepsilon)=a_0-a_1\in \Z[\beta]$.
Значит, $T_r(\varepsilon)$ построимо.
QED

\begin{pr}*\label{lower}
(a)  Если $n$ простое, то из числа 1
можно получить множество чисел, содержащее $\varepsilon_n$, используя четыре арифметические
операции и извлечения корней $(n-1)$-й степени (при которых получаются все $n-1$ значений корня).

(b) {\it Теорема Гаусса о понижении.}
Для любого $n$ из числа 1 можно получить множество чисел, содержащее $\varepsilon_n$, используя четыре арифметические операции и извлечения корней только степеней, строго меньших $n$
(при которых получаются все значения корня).
\end{pr}


\subsection{Эффективное доказательство построимости в теоремах Гаусса}\label{gaueff}


Здесь приводятся другие доказательства построимости в теореме Гаусса и теоремы Гаусса о понижении \ref{lower}.b.
Они сложнее вышеприведенных, но дают более реальную возможность получить явные формулы [BK, Sa].
Именно они принадлежит Гауссу.

Интересно бы получить явные формулы и при помощи вышеприведенного метода.

\smallskip
{\it Эффективное доказательство построимости в теореме Гаусса для $n=5$.}
Сразу выразить число $\varepsilon:=\varepsilon_5$ через радикалы трудно, поэтому сначала выразим
некоторые `многочлены от $\varepsilon$'.
Мы знаем, что $\varepsilon+\varepsilon^2+\varepsilon^3+\varepsilon^4=-1$.
Поэтому
$$(\varepsilon+\varepsilon^4)(\varepsilon^2+\varepsilon^3)=\varepsilon+\varepsilon^2+\varepsilon^3+\varepsilon^4=-1.$$
Обозначим
$$T_0:=\varepsilon+\varepsilon^4\quad\mbox{и}\quad T_1:=\varepsilon^2+\varepsilon^3.$$
Тогда по теореме Виета числа $T_0$ и $T_1$ являются корнями уравнения $t^2+t-1=0$.
Поэтому можно выразить $T_0$ (и $T_1$).
Поскольку $\varepsilon\cdot \varepsilon^4=1$, то по теореме Виета числа
$\varepsilon$ и $\varepsilon^4$ являются корнями уравнения $t^2-T_0t+1=0$.
Поэтому можно выразить $\varepsilon$ (и $\varepsilon^4$).

\smallskip
{\it Идея доказательства теоремы Гаусса о понижении \ref{lower}.b для общего случая.}
Достаточно доказать для простого $n$.
Разложим число $n-1$ в произведение $q_1q_2\dots q_s$ простых чисел.
Обозначим $\varepsilon:=\varepsilon_n$.
Сначала хорошо бы разбить сумму
$$\varepsilon+\varepsilon^2+\dots+\varepsilon^{n-1}=-1$$
на $q_1$ слагаемых $T_0,T_1,\dots,T_{q_1-1}$, которые выражаются в радикалах (иными
словами, {\it сгруппировать} хитрым образом корни уравнения
$1+x+x^2+\dots+x^{n-1}=0$).
Затем хорошо бы разбить каждую сумму $T_k$ на $q_2$ слагаемых $T_{k,0}+T_{k,1}+\dots+T_{k,q_2-1}$,
которые выражаются в радикалах.
 И так далее, пока не получим $T_{\underbrace{1,\dots,1}_s}=\varepsilon$.

Однако придумать нужные группировки чисел $1,\varepsilon,\varepsilon^2,\dots,\varepsilon^{n-1}$ нетривиально.

\begin{pr}\label{eff7}
(a) Разбейте числа $\varepsilon_7,\varepsilon_7^2,\dots,\varepsilon_7^6$ на 2 группы по 3 элемента и выразите
сумму чисел в каждой группе через квадратные радикалы.

(b)  Разбейте числа $\varepsilon_{11},\varepsilon_{11}^2,\dots,\varepsilon_{11}^{10}$ на 2 группы по 5 элементов
и выразите сумму чисел в каждой группе через квадратные радикалы.
\end{pr}

\begin{pr}\label{eff11}
Получите формулы для \quad (a) $\varepsilon_7$, \quad (b) $\varepsilon_{11}$, \quad
в которых извлекаются корни только степеней 2, 3, 5.
\end{pr}

\begin{pr}*\label{eff13}
Получите формулы для следующих чисел такие, в которых извлекаются корни только степеней 2 и 3.

(a) $\varepsilon_{13}+\varepsilon_{13}^4+\varepsilon_{13}^3+\varepsilon_{13}^{12}+\varepsilon_{13}^9+\varepsilon_{13}^{10}$.
\qquad
(b) $\varepsilon_{13}+\varepsilon_{13}^3+\varepsilon_{13}^9$.
\qquad
(c) $\varepsilon_{13}$.
\end{pr}

{\it Эффективное доказательство теоремы Гаусса о построимости.}
Достаточно доказать построимость числа $\varepsilon_n$ для $n$ простого, $n-1=2^m$.
Пусть $a_k\in\Z_2=\{0,1\}$ для каждого $k\in\{0,1,2,\dots,m-1\}$.
Обозначим
$$\overline{a_{m-1}\dots a_1a_0}:=a_0+2a_1+2^2a_2+\dots+a_{m-1}2^{m-1}.$$
Для доказательства важно, что в начале могут стоять нулевые цифры.

Обозначим через $g$ первообразный корень по модулю $n$.
Для $A\in\Z_2^k$ обозначим
$$T_A:=\sum_{B\in\Z_2^{m-k}} \varepsilon_n^{g^{\overline{BA}}}.$$
С помощью индукции по $k$ докажем, что
{\it для любых $k$ и $A\in\Z_2^k$ число $T_A$ выражается через квадратные радикалы.}
Тогда для $k=m$ получим выражаемость числа $T_{\underbrace{0\dots0}_m}=\varepsilon_n$.

База $k=0$ следует из $T_\emptyset=-1$.
Докажем шаг индукции.
Возьмем любое $A\in\Z_2^k$.
Тогда $T_A=T_{0A}+T_{1A}$.
Кроме того,
$$T_{0A}T_{1A}=
\sum_{s=0}^{n-1} N_s\varepsilon_n^s \overset{(*)}=
N_0+\sum\limits_{C\in\Z_2^k}N_{g^{\overline C}}T_C.$$
Здесь число $N_s$ (зависящее от $A$) есть количество упорядоченных решений
$$B_0,B_1\in \Z_2^{m-k-1}\quad\text{сравнения}\quad
 g^{\overline{B_00A}}+g^{\overline{B_11A}}\equiv s\mod n.$$
Ясно, что $N_s=N_{g^{2^k}s}$.
Отсюда вытекает равенство (*).
\footnote{Вот чуть более сложное окончание доказательства, которое поможет понять обобщение --- теорему Гаусса о понижении \ref{lower}.b.
Вместо $T_{0A}T_{1A}$ рассмотрим
$$(T_{0A}-T_{1A})^2=
\left(\sum_{Bl\in\Z_2^{m-k}}(-1)^l\varepsilon_n^{g^{\overline{BlA}}}\right)^2=
\sum_{l=0,s=0}^{1,n-1} |s,l|(-1)^l\varepsilon_n^s \overset{(*)}=
\sum\limits_{l=0}^1(-1)^l\left(|0,l|+\sum\limits_{C\in\Z_2^k}|g^{\overline C},l|T_C\right).$$
Здесь число $|s,l|$ (зависящее от $A$) есть количество упорядоченных решений
$$l_1,l_2\in\Z_2,\quad B_1,B_2\in \Z_2^{m-k-1}\quad\text{системы сравнений}\quad
\begin{cases}l_1+l_2 \equiv l\mod 2\\
g^{\overline{B_1l_1A}}+g^{\overline{B_2l_2A}}\equiv s\mod n\end{cases}.$$
Ясно, что $|s,l|=|g^{2^k}s,l|$.
Отсюда вытекает равенство (*).}
Значит, $T_{0A}$ и $T_{1A}$ выражаются через квадратные радикалы.
QED

\smallskip
{\it Эффективное доказательство теоремы Гаусса о понижении \ref{lower}.b.}
Пусть $n-1=q_1q_2\dots q_m$ --- разложение на простые множители (не обязательно различные).
Пусть
$$a_i\in\{0,1,2,\dots,q_{i+1}-1\} \quad\text{для каждого}\quad i\in\{0,1,2,\dots,m-1\}.$$
$$\text{Обозначим}\quad \overline{a_{m-1}\dots a_1a_0}:=a_0+a_1q_1+a_2q_1q_2+\dots+a_{m-1}q_1q_2\dots q_{m-1}.$$
(`Запись в системе счисления с переменным основанием.')
Для дальнейшего важно, что в начале могут стоять нулевые `цифры'.
Обозначим через
$$[k,l]:=\Z_{q_k}\times\dots\times\Z_{q_l}$$ множество наборов из $(k-l+1)$ `цифр', которые могут стоять в записи
$\overline{a_{m-1}\dots a_1a_0}$ на местах с $k$-го справа по $l$-е справа ($a_0$ считается первым справа).

Обозначим через $g$ первообразный корень по модулю $n$.
Для $A\in[k,1]$ обозначим
$$T_A:=\sum_{B\in[m,k+1]}\varepsilon_n^{g^{\overline{BA}}}.$$
С помощью индукции по $k$ покажем, как {\it для любых $k$ и $A\in[1,k]$ число $T_A$ выразить через радикалы.}
Тогда для $k=m$ получим выражение числа $T_{\overline{\underbrace{0\dots0}_m}}=\varepsilon_n$.

База $k=0$ следует из $T_\emptyset=-1$.
Докажем шаг индукции.
Обозначим
$q:=q_{k+1}$ и $\beta:=\varepsilon_q$.
Для любых $r=0,1,2,\dots,q-1$ и $A\in[k,1]$ обозначим
$$T_A^{(r)}:=T_{0A}+\beta^rT_{1A}+\beta^{2r}T_{2A}+\dots+\beta^{(q-1)r}T_{(q-1)A}.$$
(Резольвенты Лагранжа.)
Тогда
$$T_A^{(0)}=T_A\quad\text{и}\quad qT_{lA}=\beta^{-l}T_A^{(0)}+\beta^{-2l}T_A^{(1)}+\dots+\beta^{-(q-1)l}T_A^{(q-1)}.$$
Для любых $r=1,2,\dots,q-1$ и $A\in[k,1]$
имеем
$$(T_A^{(r)})^q=\left(\sum_{Bl\in[m,k+1]}\beta^{lr}\varepsilon_n^{g^{\overline{BlA}}}\right)^q=
\sum_{l=0,s=0}^{q-1,n-1} |s,l|\varepsilon_n^s\beta^l \overset{(*)}=
\sum_{l=0}^{q-1}\beta^l\left(|0,l|+\sum_{C\in[k,1]}|g^{\overline C},l|T_C\right).$$
Здесь число $|s,l|$ (зависящее от $A$) есть количество упорядоченных решений
$l_1,l_2,\dots,l_q\in\Z_q,\quad B_1,B_2,\dots,B_q\in[m,k+2]$ системы сравнений
$$ \begin{cases}r(l_1+l_2+\dots+l_q) \equiv l\mod q\\
g^{\overline{B_1l_1A}}+g^{\overline{B_2l_2A}}+\dots +g^{\overline{B_ql_qA}}\equiv s\mod n\end{cases}.$$
Ясно, что $|s,l|=|g^{q_1q_2\dots q_k}s,l|$.
Отсюда вытекает равенство (*).
Так $T_{lA}$ выражаются через радикалы.
QED

\subsection{Указания и решения к задачам}

 {\bf \ref{g-eps}.} Вытекает из
$$\varepsilon_n=\cos\frac{2\pi}n+\sqrt{\sin^2\frac{2\pi}n-1},\quad \cos\frac{2\pi}n=\frac{\varepsilon_n+\varepsilon_n^{-1}}2\quad\text{и}
\quad\sin\frac{2\pi}n=\frac{\varepsilon_n-\varepsilon_n^{-1}}2.$$
Или из задачи \ref{g-com}.

\smallskip
{\bf \ref{g-com}.} Часть `тогда' очевидна.
Для доказательства части `только тогда' напишите $\sqrt{a+bi}=u+vi$ и выразите $u,v$ через $a$ и $b$ с
помощью четырех арифметических операций и квадратных радикалов.

\smallskip
{\bf \ref{postr7}.} (b')
$\varepsilon_7^6+\varepsilon_7^5+\dots+\varepsilon_7+1=0$.
Как решать алгебраические уравнения $n$-й степени, у которых коэффициенты при $k$-й и при $(n-k)$-й степенях равны?

\smallskip
{\bf \ref{postr7}.bcd, \ref{7}.} Аналогично приведенному доказательству построимости числа $\varepsilon_5$.
См. подробности в следующем пункте.

\smallskip
{\bf \ref{lower}.}
(a) Аналогично доказательству построимости в теореме Гаусса.

(b) 
Докажем теорему при помощи индукции по $n$.

Если $n=ab$ для некоторых целых $0<a,b<n$, то шаг индукции следует из $\varepsilon_n=\sqrt[a]\varepsilon_b$.

Если же $n$ простое, то шаг индукции следует из (a).

\smallskip
{\bf \ref{eff7}.} (a)
Обозначим $\varepsilon:=\varepsilon_7$.
Выразим
$$T_0=\varepsilon^{3^0}+\varepsilon^{3^2}+\varepsilon^{3^4}\quad\text{и}\quad T_1=\varepsilon^3+\varepsilon^{3^3}+\varepsilon^{3^5}.$$
{\it Малоэффективное доказательство выразимости.}
Число $T_0T_1$ есть многочлен от $\varepsilon$ с целыми коэффициентами степени меньше 7
(точнее, значение в точке $\varepsilon$ некоторого многочлена от $x$ по модулю $x^7-1$ с целыми коэффициентами).
При замене $\varepsilon\to\varepsilon^3$ числа $T_0$ и $T_1$ меняются местами.
Поэтому при замене $\varepsilon\to\varepsilon^3$ число $T_0T_1$ не изменяется.
Значит, коэффициент многочлена при $\varepsilon^s$ равен его коэффициенту при $\varepsilon^{3s}$.
Так как 3 --- первообразный корень по модулю 7, то
все коэффициенты многочлена, кроме свободного члена, равны.
Из этого и $\varepsilon+\varepsilon^2+\dots+\varepsilon^6=-1$ вытекает, что $T_0T_1$ целое число.
Поэтому $T_0$ и $T_1$ выражаются.

{\it Эффективное доказательство выразимости.}
Имеем $T_0T_1=\sum_{s=0}^6N(s)\varepsilon^s$, где $N_s$ есть количество решений $(n,m)\in\Z_3^2$ сравнения
$3^{2n}+3^{2m+1}\equiv s\mod 7.$
Ясно, что $N_0+N_1+N_2+\dots+N_6=9$.
Нетрудно проверить, что $N_s=N_{3s}$.
Поэтому $N_1=N_2=\dots=N_6$.
Так как $3^0+3^1\not\equiv 0\mod7$, то $N_0\ne9$.
Из всего этого следует, что $N_0=3$ и $N_1=1$.
Значит, $T_0T_1=3-1=2$.
Отсюда $\{T_0,T_1\}=\{\dfrac{-1-\sqrt7}2,\dfrac{-1+\sqrt7}2\}$.


\smallskip
{\bf \ref{eff11}.} (a) Обозначим $\varepsilon:=\varepsilon_7$.
Имея вычисленные выше $T_0$ и $T_1$, выразим $\varepsilon$.
Обозначим
$$\beta:=\varepsilon_3\quad\text{и}\quad T_{01}:=\varepsilon^{3^0}+\beta\varepsilon^{3^2}+\beta^2\varepsilon^{3^4}.$$
{\it Малоэффективное доказательство выразимости.}
Число $T_{01}^3$ есть многочлен от $\varepsilon$ с коэффициентами в $\Z[\beta]$ степени меньше 7.
При замене $\varepsilon\to\varepsilon^{3^2}$ число $T_{01}^3$ не изменяется.
Значит, коэффициент многочлена при $\varepsilon^s$ равен его коэффициенту при $\varepsilon^{3^2s}$.
Так как 3 --- первообразный корень по модулю 7, то
коэффициенты многочлена при степенях $3^{2n}$ равны и
коэффициенты многочлена при степенях $3^{2n+1}$ равны.
Из этого и определений $T_0$ и $T_1$ вытекает, что $T_{01}^3$ выражается при помощи кубического корня
через множества $\Z[\beta,T_0,T_1]$.
Поэтому $T_{01}$ выражается как надо.

{\it Эффективное доказательство выразимости.}
Имеем $T_{01}^3=\sum_{s=1}^7\sum_{l=0}^2|s,l|\varepsilon^s\beta^l$, где $|s,l|$ есть количеству решений $(l_1,l_2,l_3)\in\Z_3^3$ системы сравнений
$$\begin{cases}l_1+l_2+l_3 \equiv l\mod 3\\
3^{2l_1}+3^{2l_2}+3^{2l_3} \equiv s\mod 7\end{cases}.$$
Ясно, что $\sum_{s=1}^7\sum_{l=0}^2|s,l|=27$.
Нетрудно проверить, что $|s,l|=|3^2s,l|$.
Поэтому $|1,l|=|2,l|=|4,l|$ и $|3,l|=|5,l|=|6,l|$.
Значит, $T_{01}^3=A+BT_0+CT_1$ для некоторых  $A, B, C\in\Z[\beta]$.
Их несложно найти.

{\it Окончание указаний.}
Аналогично выражается $T_{02}:=\varepsilon^{3^1}+\beta\varepsilon^{3^3}+\beta^2\varepsilon^{3^5}$.
Потом выражается $\varepsilon=\dfrac{T_0+T_{01}+T_{02}}3$.

\smallskip
{\bf \ref{eff13}.} (a) Обозначим $\varepsilon=\varepsilon_{13}$,
$$A_0:=\varepsilon+\varepsilon^4+\varepsilon^3+\varepsilon^{12}+\varepsilon^9+\varepsilon^{10}\quad\text{и}
\quad A_1:=\varepsilon^2+\varepsilon^5+\varepsilon^6+\varepsilon^7+\varepsilon^8+\varepsilon^{11}.$$
Тогда $A_0+A_1=-1$ и $A_0A_1=-3$.
Значит, $\{A_0,A_1\}=\{\dfrac{-1-\sqrt{13}}2,\dfrac{-1+\sqrt{13}}2\}.$
Из рисунка правильного 13-угольника (или из оценок) видно, что  $A_1<0$.
Значит, $A_1=\dfrac{-1-\sqrt{13}}2$ и $A_0=\dfrac{-1+\sqrt{13}}2$.


\comment

\begin{pr}*\label{expl}
(a) Получите явную формулу для $\varepsilon_{11}$, использующую целые числа, четыре арифметические операции и извлечения корней второй и пятой степени.
Важным дополнением к формуле будет указание, какие именно значения корней берутся.

(b) Найдите наименьшее количество корней пятой степени, которое может быть в формуле из (a).
\end{pr}

\bigskip
{\bf Приложение: другое доказательство теоремы Гаусса о понижении \ref{lower}.b}

\smallskip
{\it  Набросок другого доказательства построимости в теореме Гаусса.}
Начало до формулы $T_A=T_{\overline{0A}}+T_{\overline{1A}}$ такое же, как в эффективном доказательстве.
При замене $x\to x^{g^{2^k}}$ `многочлены'  $T_{\overline{0A}}$ и $T_{\overline{1A}}$ по модулю $x^n-1$ меняются местами, поэтому многочлен $T_{\overline{0A}}T_{\overline{1A}}$ по модулю $x^n-1$ не меняется.
Значит, коэффициент при $x^s$ соответствующего многочлена степени меньше $n$ равен его коэффициенту при $x^{sg^{2^k}}$ (степени рассматриваются по модулю $n$).
Так как $g$ --- первообразный корень по модулю $n$, то для любого фиксированного $A'\in\Z_2^k$ и различных
$B\in\Z_2^{m-k}$ коэффициенты этого многочлена при $x^{g^{\overline{BA'}}}$ равны.
Из этого и определения многочленов $T_A$ вытекает, что число $T_{\overline{0A}}(\varepsilon_n)T_{\overline{1A}}(\varepsilon_n)$ является линейной комбинацией с целыми коэффициентами чисел 1 и $T_{A'}(\varepsilon_n)$ с $A'\in[k,1]$.
Это влечет шаг индукции.   QED



\smallskip
{\it Набросок другого доказательства теоремы Гаусса о понижении \ref{lower}.b.}
Начало такое же, как в эффективном доказательстве.
Число $T_{A,r}^q$ есть многочлен от $\alpha$ с коэффициентами в $\Z[\varepsilon]$ степени меньше $p$.
(Точнее, значение в точке $\alpha$ некоторого многочлена с коэффициентами в $\Z[\varepsilon]$ от $x$ по модулю $x^p$.)
При замене $\alpha\to\alpha^{g^{q_1q_2\dots q_k}}$ это число $T_{A,r}^q$ не изменяется.
Значит, коэффициент этого многочлена при $\alpha^s$ равен его коэффициенту при $\alpha^{sg^{q_1q_2\dots q_k}}$.
Так как $g$ --- первообразный корень по модулю $p$, то для любого фиксированного $A'\in[k,1]$ и различных
$B\in[m,k+1]$ коэффициенты многочлена при $\alpha^{g^{\overline{BA'}}}$ равны.
Из этого и определения чисел $T_A$ вытекает, что
число $T_{A,r}^q$ является линейной комбинацией с коэффициентами в $\Z[\varepsilon]$ чисел 1 и $T_{A'}$ с $A'\in[k,1]$.
Шаг индукции следует из этого и формул
$$qT_{\overline{lA}}=
T_{A,0}+\varepsilon^{-l}T_{A,1}+\varepsilon^{-2l}T_{A,2}+\dots+\varepsilon^{-(q-1)l}T_{A,q-1}.\quad QED$$

\begin{pr}\label{glemma}
Пусть $p,q$ простые, $g$ первообразный корень по модулю $p$ и $qd$ делитель числа $p-1$.
Обозначим
$$\alpha:=\varepsilon_p
\quad\text{и}\quad t_{d,u}:=\alpha^{g^u}+\alpha^{g^{d+u}}+\alpha^{g^{2d+u}}+\dots+\alpha^{g^{p-1-d+u}}\quad\text{для}\quad u=0,1,2,\dots,d-1.$$
Тогда существуют $a_1,\dots,a_s\in K_d:=\Q[t_{d,1},t_{d,2},\dots,t_{d,d-1}]$, для которых
\linebreak
$t_{qd,v}\in K_d[\varepsilon_q,\sqrt[q]{a_1},\dots,\sqrt[q]{a_s}]$ при любом $v=0,1,2,\dots,qd-1$
(т.е. $t_{qd,v}$ получается на калькуляторе с использованием $t_{d,1},t_{d,2},\dots,t_{d,d-1}$ и нескольких извлечений корня $q$-й степени).
\end{pr}

\endcomment

\section{Задачи о неразрешимости в радикалах}\label{non}


В параграфах 4 и 5 через $\Q$ обозначается множество всех рациональных чисел;
`многочлен с рациональными коэффициентами' коротко называется многочленом.
Многочлен называется {\it неприводимым} над множеством $F$, если он не раскладывается в произведение многочленов меньшей степени
с коэффициентами в $F$.

\subsection{Одно извлечение квадратного корня}\label{non12}

\begin{pr}\label{number2} Представимо ли следующее число в виде $a+\sqrt b$, где $a,b\in\Q$?

(a) $\sqrt{3+2\sqrt2}$; \quad (a') $\sqrt{2+\sqrt2}$; \quad (b) $\sqrt[3]{\sqrt 5+2}-\sqrt[3]{\sqrt 5-2}$;
\quad (c) $\sqrt[3]{7+5\sqrt2}$; \quad
\linebreak
(d) $\cos(2\pi/5)$; \quad (e) $\sqrt[3]2$; \quad (f) $\sqrt2+\sqrt[3]2$; \quad (g) $\cos(2\pi/9)$; \quad (h) $\cos(2\pi/7)$.
\end{pr}

\begin{pr}\label{3alpha} Число $\cos(2\pi/9)$ является корнем уравнения $8x^3-6x+1=0$.
\end{pr}

\begin{pr}\label{sopr}
(a) {\bf Лемма о сопряжении.} Если $a,b\in\Q$ и $a+b\sqrt2$ --- корень многочлена, то $a-b\sqrt2$ --- тоже его корень.

(b) {\bf Лемма о линейной независимости.} Если $a+b\sqrt2=0$ для некоторых $a,b\in\Q$, то $a=b=0$.

(c) {\bf Утверждение.} Если многочлен степени выше второй неприводим над $\Q$,
 то ни один из его корней не представим в виде $a\pm\sqrt b$, где $a,b\in\Q$.
\end{pr}

\begin{pr}\label{calc2}
{\bf Лемма о калькуляторе.} Пусть $F\in\{\R,\C\}$.
Число, которое можно получить  на $F$-калькуляторе так, чтобы извлечение корня происходило только один раз, причем второй степени,
имеет вид $a\pm\sqrt b$, где $a,b\in\Q$ и, для $F=\R$, $b>0$.
\end{pr}

Задачи \ref{number2} и \ref{sopr}.c интересны в связи с неразрешимостью в радикалах, поскольку нам нужно придумать многочлен,
корни которого невозможно получить на калькуляторе, а числа из задачи \ref{number2} являются корнями многочленов (подумайте, каких).

Из утверждения \ref{sopr}.c и леммы \ref{calc2} о калькуляторе вытекает, что
{\it если многочлен степени выше второй неприводим над $\Q$, то ни один из его корней невозможно получить на вещественном калькуляторе так, чтобы извлечение корня происходило только один раз, причем второй степени}.
Это --- наше первое продвижение к теоремам о неразрешимости в радикалах.
Аналогичные продвижения в следующих двух пунктах (сформулируйте их самостоятельно) вытекают из соответствующих
утверждений и лемм о калькуляторе.

\begin{pr}\label{quairr} Для каких $n$ число $\cos(2\pi/n)$

(a) рационально? \qquad (b)* представимо в виде $a+\sqrt b$, где $a,b\in\Q$?
\end{pr}

\subsection{Одно извлечение корня четвертой степени}\label{non14}

\begin{pr}\label{foex} Представимо ли следующее число в виде $a+b\sqrt[4]2+c\sqrt2+d\sqrt[4]8$, где  $a,b,c,d\in\Q$?

(a) $\sqrt[3]3$; \quad (b) $\sqrt[6]3$; \quad (c) $\sqrt[4]3$.
\end{pr}

\begin{pr}\label{sopr4}
{\bf Лемма о сопряжении.} Пусть $a,b,c,d\in\Q$ и $a+b\sqrt[4]2+c\sqrt2+d\sqrt[4]8$ корень многочлена.
Тогда корнем этого многочлена также является число

(a) $a-b\sqrt[4]2+c\sqrt2-d\sqrt[4]8$; \qquad
(b) $a-c\sqrt2+i\sqrt[4]2(b-d\sqrt2)$ и $a-c\sqrt2-i\sqrt[4]2(b-d\sqrt2)$.
\end{pr}

\begin{pr}\label{sopr4det} {\bf Лемма о линейной независимости.}

(a)  Если $a+b\sqrt[4]2+c\sqrt2+d\sqrt[4]8=0$ для некоторых $a,b,c,d\in\Q$, то $a=b=c=d=0$.

(b) Если $a+b\sqrt[4]2+c\sqrt2+d\sqrt[4]8=0$ для некоторых $a,b,c,d\in\Q[i]:=\{x+iy\ :\ x,y\in\Q\}$, то $a=b=c=d=0$.
\end{pr}

\begin{pr}\label{utv4} {\bf Утверждение.} Если многочлен степени, отличной от 1, 2 и 4, неприводим над $\Q$, то
ни один из его корней не представим в виде $a+br+cr^2+dr^3$, где $r\in\C$ и $a,b,c,d,r^4\in\Q$.
 \end{pr}

\begin{pr}\label{mno4}
{\bf Лемма о калькуляторе.} Пусть $F\in\{\R,\C\}$. Число, которое можно получить  на $F$-калькуляторе так, чтобы извлечение корня
происходило только один раз, причем четвертой степени, имеет вид $a+br+cr^2+dr^3$, где $r\in F$ и $a,b,c,d,r^4\in\Q$.
\end{pr}

\subsection{Несколько извлечений квадратных корней }\label{nonn2}

\begin{pr}\label{twoexa} Существуют ли рациональные числа $a,b,c,d$, для которых $\sqrt[3]2=$

(a) $\displaystyle\frac{a+\sqrt b}{c+\sqrt b}$;
\quad (b) $a+\sqrt b+\sqrt c$; \quad
(c) $a+\sqrt{b+\sqrt c}$; \quad (d) $a+\sqrt b+\sqrt c+\sqrt d$?
\end{pr}


\begin{pr}\label{noncon}
Следующие числа не являются вещественно построимыми:
\linebreak
(a) $\sqrt[3]2$; \quad (b) $\cos(2\pi/9)$; \quad
(c) произвольный корень уравнения $y^3+y+1=0$;
\quad
(d) $\cos(2\pi/7)$.
\end{pr}

 \begin{pr}\label{tower2}
(a) Оторвем у комплексного калькулятора кнопку '$:$', но разрешим использовать все рациональные числа.
Тогда множество чисел, которые можно получить на калькуляторе, не изменится.

(b) {\bf Лемма о калькуляторе.} Для $F\subset\C$, \ $r\in\C$ и $r^2\in F$ обозначим
\linebreak
$F[r]:=\{a_0+a_1r\ |\ a_0,a_1\in F\}$.

Число $x$ построимо тогда и только тогда, когда существуют такие $r_1,\ldots r_{s-1}\in\C$, что
$$\Q=F_1\subset F_2\subset F_3\subset \ldots\subset F_{s-1}\subset F_s\ni x,
\quad\mbox{где}\quad r_k^2\in F_k,\quad r_k\not\in F_k\quad\mbox{и}\quad F_{k+1}=F_k[r_k]$$
для любого $k=1,\dots,s-1$.

(Такая последовательность называется {\it башней квадратичных расширений}.
Доказательство невозможности, основанное на рассмотрении аналогичных цепочек-башен,
называется в математической логике и программировании {\it индукцией по глубине формулы}.)

(c)
Некоторый (или, эквивалентно, каждый) корень кубического многочлена построим тогда и только тогда,
когда один из корней этого многочлена рационален.

(d) Некоторый (или, эквивалентно, каждый) корень многочлена
4-й степени построим тогда и только тогда,
когда его {\it кубическая резольвента} (\S\ref{mot34})
имеет рациональный корень.
\end{pr}

\begin{pr}\label{alge}
Любое построимое число (и даже любое число, которое можно получит на комплексном калькуляторе)
является {\it алгебраическим}, т.е. корнем некоторого многочлена с целыми коэффициентами.
(Из этого и доказанной в 1883 г. Линдеманом трансцендентности числа $\pi$, влекущей трансцендентность числа $\sqrt{\pi}$, вытекает, что задача о квадратуре круга неразрешима циркулем и линейкой.)
\end{pr}

\subsection{К доказательству непостроимости в теореме Гаусса}\label{nongau}

\begin{pr}\label{irre} Найдите неприводимый над $\Q$ многочлен, корнем которого является число

(a) $\sqrt2+\sqrt3$; \quad (b) $\sqrt2+\sqrt3+\sqrt5$; \quad
(c) $\sqrt{2+\sqrt3}$; \quad
(d) $\sqrt{5+\sqrt6}$; \quad (e) $\sqrt{5+\sqrt{1+\sqrt3}}$. \quad
\end{pr}

\begin{pr}\label{baslem}
(a) Если неприводимый над $\Q[\sqrt2]$ многочлен $P$ имеет корень вида $a+\sqrt b$, где $a,b\in \Q[\sqrt2]$, то
$\deg P\in\{1,2\}$.

(b) Если неприводимый над $\Q$ многочлен $P$ имеет корень вида $a+\sqrt b+\sqrt c$, где $a,b,c\in \Q$, то
$\deg P\in\{1,2,4\}$.

(c) Если неприводимый над $\Q$ многочлен $P$ имеет корень вида $\sqrt a+\sqrt{b+\sqrt c}$, где $a,b,c\in \Q$, то
$\deg P\in\{1,2,4,8\}$.

(d) Если неприводимый над $\Q$ многочлен $P$ имеет построимый корень, то $\deg P$ есть степень двойки.
\end{pr}

\begin{pr}\label{irree} Найдите неприводимый над $\Q$ многочлен, корнем которого является число

(5)  $\varepsilon_5$; \quad (7) $\varepsilon_7$; \quad (9)  $\varepsilon_9$; \quad (11) $\varepsilon_{11}$; \quad
(13) $\varepsilon_{13}$; \quad (25) $\varepsilon_{25}$.
\end{pr}

\begin{pr}
(a) {\it Лемма Гаусса.} Если многочлен с целыми коэффициентами неприводим над $\Z$, то он неприводим и над $\Q$.

(b) {\it Признак Эйзенштейна.} Пусть $p$ простое. Если для многочлена с целыми коэффициентами старший коэффициент  не делится на $p$, остальные делятся на $p$, а свободный член не делится на $p^2$, то этот многочлен неприводим над $\Z$.
\end{pr}

\begin{pr}\label{incons}
Докажите непостроимость в теореме Гаусса.
\end{pr}

 \subsection{Одно извлечение корня третьей степени}\label{non13}

\begin{pr}\label{3ex} Представимо ли следующее число в виде $a+b\sqrt[3]2+c\sqrt[3]4$, где $a,b,c\in\Q$?

(a) $\sqrt3$; \quad (b) $\cos(2\pi/9)$;  \quad (c) $\sqrt[5]3$; \quad (d) $\sqrt[3]3$.

(e) наименьший положительный корень уравнения $x^3-4x+2=0$.
\end{pr}


\begin{pr}\label{sopr3}
Пусть $\varepsilon:=\varepsilon_3=\cos\dfrac{2\pi}3+i\sin\dfrac{2\pi}3$,\ $r\in\R-\Q$ и $r^3,a,b,c\in\Q$.

(a) {\bf Лемма о сопряжении.} Если многочлен имеет корень $a+br+cr^2$, то корнями этого многочлена являются также числа
$a+b\varepsilon r+c\varepsilon^2r^2$ и $a+b\varepsilon^2r+c\varepsilon r^2$.


(b) Если $a+br+cr^2=0$, то $a=b=c=0$.

(c) {\bf Лемма о линейной независимости.} Если $k+lr+mr^2=0$ для некоторых
$k,l,m\in\Q[\varepsilon]:=\{x+y\varepsilon\ :\ x,y\in\Q\}$, то $k=l=m=0$.


(d) {\bf Лемма о рациональности.}
Кубический многочлен с корнями
$$(*)\quad a+br+cr^2, \quad a+br\varepsilon+cr^2\varepsilon^2\quad\text{и}\quad a+br\varepsilon^2+cr^2\varepsilon$$
и коэффициентом 1 при $x^3$ имеет рациональные коэффициенты.
\end{pr}

\begin{pr}\label{utv3}
{\bf Утверждение.} Если многочлен неприводим над $\Q$ и имеет корень вида $a+br+cr^2$, где $r\in\R-\Q$ и $a,b,c,r^3\in\Q$, то
степень многочлена равна 3 и он имеет ровно один вещественный корень.
\end{pr}


\begin{pr}\label{que3} (3R) Как по многочлену третьей степени узнать, имеет ли он корень вида $a+br+cr^2$, где $a,b,c,r^3\in\Q$ и $r\in\R$?

(3C) То же для $r\in\C$.

(3CC) То же для $a,b,c,r^3\in\Q[i]$ и $r\in\C$?

(3Rn), (3Cn), Те же вопросы для корня вида
$a_0+a_1r+a_2r^2+\dots+a_{n-1}r^{n-1}$, где $n\in\Z$, $n\ge1$, $a_0,a_1,\dots,a_{n-1},r^n\in\Q$ и $r\in\R$?
(Иными словами, как по многочлену третьей степени узнать, имеет ли он корень,
который можно получить на вещественном/комплексном  калькуляторе так, чтобы извлечение корня происходило только один раз?)

(4R)*, (4C)*, (4Rn)*, (4Cn)* Те же вопросы для многочленов четвертой степени.
\end{pr}


\begin{pr}\label{mno3} {\bf Лемма о калькуляторе.} Пусть $F\in\{\R,\C\}$.
Число, которое можно получить  на $F$-калькуляторе так, чтобы извлечение корня происходило только один раз,
причем третьей степени, имеет вид $a+br+cr^2$, где $r\in F$ и $a,b,c,r^3\in\Q$.
 \end{pr}

 \subsection{Одно извлечение корня простой степени}\label{non1p}

\begin{pr}\label{kvadiffreal} Представимо ли следующее число в виде $a_0+a_1\sqrt[7]2+a_2\sqrt[7]{2^2}+\dots+a_6\sqrt[7]{2^6}$, где
$a_0,a_1,a_2,\dots,a_6\in\Q$?

(a) $\sqrt3$; \qquad (b) $\cos(2\pi/21)$; \qquad (с) $\sqrt[11]3$; \qquad (d) $\sqrt[7]3$;

(e) наименьший положительный корень уравнения $x^7-4x+2$.
\end{pr}

\begin{pr}\label{idea} Пусть $q$ нечетное простое, $\varepsilon_q:=\cos\dfrac{2\pi}q+i\sin\dfrac{2\pi}q$,
\ $A$ --- многочлен степени меньше $q$, \ $r\in\R-\Q$ и $r^q\in\Q$.

(a) {\bf Слабая лемма о неприводимости.} Многочлен $x^q-r^q$ неприводим над $\Q$.

(b) {\bf Слабая лемма о линейной независимости.} Если $A(r)=0$, то $A=0$.

(c) {\bf Лемма о сопряжении.} Если многочлен имеет корень $A(r)$, то
он имеет также корни $A(r\varepsilon_q^k)$ для каждого $k=1,2,3,\dots,q-1$.

(d) {\bf Лемма о неприводимости.} Многочлен $x^q-r^q$ неприводим над $\Q[\varepsilon_q]$.

(e) {\bf Лемма о линейной независимости.} Если $B\in\Q[\varepsilon_q][x]$ --- многочлен степени меньше $q$ с коэффициентами в $\Q[\varepsilon_q]$ и $B(r)=0$, то $B=0$.


(f) {\bf Лемма о рациональности.}
Многочлен
$(x-A(r))(x-A(r\varepsilon_q))\dots(x-A(r\varepsilon_q^{q-1}))$
имеет рациональные коэффициенты.
\end{pr}

Следующее утверждение интересно и нетривиально даже для многочленов третьей степени.

\begin{pr}\label{kronint1} {\bf Утверждение.}
Если многочлен
неприводим над $\Q$ и имеет более одного вещественного корня,
то ни один из его корней не представим в виде $A(r)$
ни для каких многочлена $A$, простого нечетного $q$ и

(a) $r\in\R$, причем $r^q\in\Q$.
\quad
(b) $r\in\C$, причем $r^q\in\Q$ и $r^q\ne b^q$ ни для какого $b\in\Q$.
\end{pr}

\begin{pr}\label{mnoq}
{\bf Лемма о калькуляторе.} Пусть $K\in\{\R,\C\}$.
Число, которое можно получить на $K$-калькуляторе так, чтобы извлечение корня происходило только один раз, равно
$A(r)$ для некоторых $r\in K$, $q\in\Z$ и $A\in\Q[x]$, причем $r^q\in\Q$.
\end{pr}

Из утверждения \ref{kronint1} и леммы \ref{mnoq} о калькуляторе вытекает, что
{\it если многочлен неприводим над $\Q$ и имеет более одного вещественного корня,
то ни один из его корней невозможно получить на вещественном калькуляторе так, чтобы извлечение корня происходило только один раз}.
Ср. с сильной веществееной теоремой о неразрешимости из \S\ref{prostr}.

 \subsection{Несколько извлечений корней }\label{nonn}

\begin{pr}\label{mailemtow}
(a-f) Докажите аналоги утверждений задачи \ref{idea} с заменой $\Q$ на произвольное подмножество $F\subset\R$,
замкнутое относительно операций сложения, вычитания, умножения и деления на ненулевое число
(и многочленов с коэффициентами в $\Q$ на многочлены с коэффициентами в $F$).
\end{pr}

\begin{pr}\label{twocalc}
(a) Можно ли число $\cos(2\pi/9)$ получить на вещественном калькуляторе так, чтобы извлечение корня происходило
только два раза, причем оба раза простой степени?

(b) Придумайте многочлен пятой степени, ни один из корней которого невозможно получить на комплексном калькуляторе
так, чтобы извлечение корня происходило только два раза, причем оба раза простой степени.
\end{pr}

Если следующие задачи
не получаются, к их решению можно вернуться после \S\ref{pro}.

\begin{pr}\label{isslreal}
(a) Один корень (или, эквивалентно, каждый корень) кубического многочлена можно получить на вещественном
калькуляторе тогда и только тогда, когда этот многочлен имеет либо хотя бы один рациональный корень, либо ровно
один вещественный корень.

(b)* {\it Гипотеза.} Каждый корень неприводимого многочлена четвертой степени можно получить на
вещественном калькуляторе тогда и только тогда, когда хотя бы один корень его кубической резольвенты
можно получить на вещественном калькуляторе.
\end{pr}

\begin{pr}\label{quairrgen}* Для каких $n$ число $\cos(2\pi/n)$ можно получить на вещественном калькуляторе?
\end{pr}

\begin{pr}\label{issl4}* Пусть $F\in\{\R,\C\}$. Как по многочлену четвертой степени узнать, имеет ли он

(2R), (2C) корень, который можно получить на $F$-калькуляторе, но чтобы извлечение корня происходило только два
раза, причем оба раза простой степени?

(3R), (3C) Те же вопросы, но чтобы извлечение корня происходило только 3,4,$\dots$ раз.
\end{pr}

\begin{pr}\label{root}
(a) Если многочлен простой степени $p$ неприводим над $\Q[\varepsilon_5]$ и приводим над $\Q[\varepsilon_5,\sqrt[5]2]$, то $p=5$ и многочлен имеет корень в $\Q[\varepsilon_5,\sqrt[5]2]$.

(b) Если многочлен простой степени $p$ неприводим над $\Q$ и приводим над $\Q[\sqrt[5]2]$,
то $p=5$ и многочлен имеет корень в $\Q[\varepsilon_5,\sqrt[5]2]$.

(c) Существуют ли целое $a$, простое $q$ и многочлен простой степени $p$, неприводимый над $\Q$ и приводимый над $\Q[\sqrt[q]a]$, не имеющий корня в $\Q[\sqrt[q]a]$?
\end{pr}

\comment

Обозначим для $F\in\{\R,\C\}$ \quad
$$\widehat F:=\{a+br+cr^2\ |\ r\in F,\ a,b,c,r^3\in\Q\}.$$
\begin{pr}\label{kvadiff} (a,b,c) Какие из чисел задачи \ref{3ex} лежат в $\widehat\R$?
\end{pr}

\begin{pr}\label{kvadiffc} (a,b,c) Какие из чисел задачи \ref{3ex} лежат в $\widehat\C$?

(d) Лежит ли $i\sqrt3$ в $\widehat\C$?
\end{pr}

(a) Квадратный многочлен имеет корень из $\widehat\C$ тогда и только тогда, когда он имеет рациональный корень.

(b) Квадратный многочлен имеет корень из $\widehat\C$ тогда и только тогда, когда
либо он имеет рациональный корень, либо его дискриминант имеет вид $-3h^2$ для некоторого $h\in\Q$.

  лежит в $\widehat F$.

\begin{pr}\label{que3i}*
\end{pr}



\smallskip
{\bf \ref{kvadiff}.}
Не лежит ни одно.
Доказательство аналогично задаче \ref{3ex}.
Приведем детали для первого решения пункта (a).

(a) Пусть можно, т.е. $\sqrt3=a+br+cr^2$ для некоторых $r\in\R$, $r\not\in\Q$, $a,b,c,r^3\in\Q$.
Тогда
$$2=(a^2+2bcr^3)+(2ab+c^2r^3)r+(2ac+b^2)r^2.$$
Так как многочлен $x^3-r^3$ не имеет рациональных корней, то он неприводим над $\Q$.
Значит, $2ab+c^2r^3=2ac+b^2=0$ (ср. \ref{sopr3}.b).
Поэтому $b^3=-2abc=c^3r^3$.
Тогда либо $b=c=0$, либо $r=b/c$.
Оба случая невозможны.

\smallskip
{\bf \ref{kvadiffc}.}
Лежит только $i\sqrt3$.
Оно получается после извлечения корня кубического из 1.

(a,b,c) Если корень извлекался комплексный и из куба рационального числа, то
\linebreak
$\sqrt3,\sqrt[5]3,\cos(\pi/9)\in\Q[\varepsilon]\cap\R=\Q$, что неверно.
Если же корень извлекался комплексный и не из куба рационального числа, то
аналогично вещественному случаю.
Для комплексного случая в леммах нужно брать $r\in\C$, но $r^3\ne p^3$ ни для какого $p\in\Q$.

\begin{pr}\label{irred}
Следующие многочлены неприводимы над $\Q$:

(a) $x^5-4$; \quad
(b) $x^q-a$, где $q$ простое, $a\in\Q$ и $\sqrt[q]a\not\in\Q$.
\end{pr}






\endcomment

\section{Указания и решения к задачам из \S\ref{non}}\label{hints}

 \subsection*{Одно извлечение квадратного корня}

 {\bf \ref{number2}.} (a,b,c,d) можно, (a',e,f,g,h) нельзя.

(a,c) $\sqrt{3+2\sqrt2}=\sqrt[3]{7+5\sqrt2}=1+\sqrt2$.

(b) $\sqrt[3]{\sqrt 5+2}-\sqrt[3]{\sqrt 5-2}=1$.

(d) $\cos(2\pi/5)=(\sqrt5-1)/4$.

(e) Пусть можно.
Тогда $2=(\sqrt[3]2)^3=(a^3+3ab)+(3a^2+b)\sqrt b$.
Так как $3a^2+b\ne0$, то $\sqrt b\in\Q$.
Значит, $\sqrt[3]2\in\Q$ --- противоречие.

Другие способы --- аналогично пунктам (e,h) или утверждению \ref{sopr}.c.


(f) Пусть можно, т.е. $\sqrt2+\sqrt[3]2=a+\sqrt b$.
Число $\sqrt2+\sqrt[3]2$ является корнем многочлена
$((x-\sqrt2)^3-2)((x+\sqrt2)^3-2)$ с рациональными коэффициентами.
Тогда по лемме о сопряжении (\ref{sopr}.a) этот многочлен имеет корень $a-\sqrt b$.
По теореме о рациональных корнях у этого многочлена нет рациональных корней.
Значит, $b\ne0$ и корни $a\pm\sqrt b$ различны.
Но у этого многочлена только два вещественных корня: $\sqrt2+\sqrt[3]2$ и $-\sqrt2+\sqrt[3]2$.
Поэтому $a+\sqrt b=\sqrt2+\sqrt[3]2$ и $a-\sqrt b=-\sqrt2+\sqrt[3]2$.
Отсюда $\sqrt[3]2=a\in\Q$. Противоречие.

Другой способ --- аналогично (e) и \ref{foex}.b.

(g) Из задачи \ref{3alpha} aналогично решению пункта (f)  получаем, что  числа $a+\sqrt b$ и $a-\sqrt b$
являются различными корнями многочлена $8x^3-6x+1$.
Тогда по теореме Виета третий корень равен $-2a\in\Q$.
Противоречие.

Другой способ --- аналогично утверждению \ref{sopr}.c.

(h) (И. Брауде-Золотарев) Из равенства $1+\varepsilon_7+\varepsilon_7^2+\dots+\varepsilon_7^6=0$ получаем
$\cos(2\pi/7)+\cos(4\pi/7)+\cos(6\pi/7)=-1/2$.
Используя формулы для $\cos2\alpha$ и $\cos3\alpha$, получаем, что число $\cos(2\pi/7)$ является корнем уравнения
$8t^3+4t^2-4t-1=0$.
Заменим $u=2t$, получим $u^3+u^2-2u-1=0$.
Это уравнение не имеет рациональных корней.
Значит, уравнение $8t^3+4t^2-4t-1=0$ тоже.
Поэтому мноочлен $8t^3+4t^2-4t-1=0$ неприводим над $\Q$.
Значит, по задаче \ref{sopr} получаем утверждение задачи.

{\it Указание к другому решению.}
Докажите, что $\varepsilon_7$ невозможно получить на калькуляторе так, чтобы корень извлекался только квадратный и только 3 раза.
Ср. с задачами \ref{twoexa} и \ref{noncon}.d.


\smallskip
{\bf \ref{3alpha}.} Выразите $\cos3\alpha$ через $\cos\alpha$.

\smallskip
{\bf \ref{sopr}.}
(a) Ввиду возможности делить с остатком на $(x-a)^2-2b^2$ утверждение достаточно доказать для многочленов первой степени.
Ввиду иррациональности числа $\sqrt2$ если $a+b\sqrt2$ является корнем многочлена первой степени, то этот многочлен нулевой.
Значит, число $a-b\sqrt2$ также является его корнем.
\footnote{Следующую интерпретацию приведенного решения (и его обобщений) предложил Л. Шабанов.
Подставим в многочлен $x=a+br$ и раскроем скобки, заменяя всюду $r^2$ на 2.
Получим выражение $m+nr$.
Подставляя $r=\sqrt2$, получаем $m+n\sqrt2=0$.
По лемме о линейной независимости (b) $m=n=0$.
Подставляя $r=-\sqrt2$, видим, что значение многочлена в точке $a-b\sqrt2$ равно $m-n\sqrt2=0$.}

(c) Если $\sqrt b\in\Q$, то утверждение очевидно.
Пусть $\sqrt b\not\in\Q$.
Поделим данный многочлен с остатком на $(x-a)^2-2b^2$.
В остатке получится многочлен первой степени, имеющий корень $a+\sqrt b$.
Так как $\sqrt b\not\in\Q$, то остаток нулевой.
Поэтому данный многочлен делится на $(x-a)^2-b$.
Противоречие с его неприводимостью над $\Q$.



\smallskip
{\bf \ref{calc2}.} Достаточно доказать, что множество чисел такого вида замкнуто
относительно сложения, вычитания, умножения, и деления.
Это, естественно, не так.

Поэтому обозначим через $\sqrt c$ число, полученное при единственном извлечении корня, где $c\in\Q$.
И будем доказывать, что тогда все полученные числа имеют вид $a+b\sqrt c$, где $a,b\in\Q$.
Достаточно доказать, что множество чисел такого вида замкнуто
относительно сложения, вычитания, умножения, и деления.
Это неочевидно только для деления, для чего оно следует из $(a+b\sqrt c)(a-b\sqrt c)=a^2-b^2c$.

\smallskip
{\bf \ref{quairr}.} (a) Ответ: $n\in\{1,2,3,4\}$.

  \subsection*{Одно извлечение корня четвертой степени}

  {\bf \ref{foex}.}  Нельзя.

(a) {\it Первое решение.} Пусть можно.
По леммам о сопряжении и о линейной независимости (\ref{sopr4}.a и \ref{sopr4det}.a)
многочлен $x^3-3$ имеет два различных вещественных корня. Противоречие.

{\it Второе решение.} Пусть можно.
По леммам о сопряжении и о линейной независимости (\ref{sopr4}.b и \ref{sopr4det}.b)
многочлен $x^3-3$ имеет четыре различных корня. Противоречие.

(b) {\it Первое решение.} Пусть можно.
По леммам о сопряжении и о линейной независимости (\ref{sopr4}.a и \ref{sopr4det}.a)
многочлен $x^6-3$ имеет два различных вещественных корня
$$a+b\sqrt[4]2+c\sqrt2+d\sqrt[4]8=\sqrt[6]3\quad\text{и}\quad a-b\sqrt[4]2+c\sqrt2-d\sqrt[4]8=-\sqrt[6]3.$$
Тогда $a=c=0$ и $b+d\sqrt2=\sqrt[6]3/\sqrt[4]2$.
Последнее число есть корень многочлена $x^{12}-9/8$.
Значит, $b-d\sqrt2$ тоже его корень.
Поэтому $b-d\sqrt2=-\sqrt[6]3/\sqrt[4]2$. Отсюда $b=0$.
Противоречие с рациональностью $d$.

{\it Второе решение.} Пусть можно.
По лемме о сопряжении (\ref{sopr4}.b) многочлен $x^6-3$ имеет четыре различных корня $x_1,x_2,x_3,x_4$,
указанные в \ref{sopr4}.b.
Так как ни один из них не рационален, то $b=c=d=0$ невозможно.
Значит, по лемме о линейной независимости (\ref{sopr4det}.b) эти корни различны.
Поэтому многочлен $x^6-3$ делится на $(x-x_1)(x-x_2)(x-x_3)(x-x_4)$.
У многочлена $(x-x_1)(x-x_2)(x-x_3)(x-x_4)$ рациональные коэффициенты (докажите!).
Противоречие с неприводимостью многочлена $x^6-3$ над $\Q$.

(c) Аналогично первому решению пункта (b).

\smallskip
{\bf \ref{sopr4}.} (a) Подставим в многочлен $x=a+bt+ct^2+dt^3$ и поделим с остатком на $t^4-2$.
Подставляя $t=\sqrt[4]2$, получаем по лемме о линейной независимости (\ref{sopr4det}.b), что остаток нулевой.
Значит, если $r^4=2$, то  $a+br+cr^2+dr^3$ есть корень исходного многочлена.



\smallskip
{\bf \ref{sopr4det}.} (a) {\it Первое решение.} Перепишем условие в виде $(a+c\sqrt2)+(b+d\sqrt2)\sqrt[4]2=0$.
Так как $b+d\sqrt2\ne0$, то $-\sqrt[4]2=\dfrac{a+c\sqrt2}{b+d\sqrt2}=A+B\sqrt2$ для некоторых $A,B\in\Q$.
Возводя в квадрат, получаем $A^2+2B^2=0$. Противоречие.

{\it Второе решение.} Так как многочлен $x^4-2$ неприводим над $\Q$, то он не может иметь общий корень с многочленом $a+bx+cx^2+dx^3$ третьей степени.

(b) Докажите отдельно для вещественной и мнимой части.







\smallskip
{\bf \ref{utv4}.} Аналогично задаче \ref{foex}.
Отдельно рассматривается более простой случай $r^2\in\Q$.

\smallskip
{\bf \ref{mno4}.} Достаточно доказать, что число, обратное к ненулевому числу такого вида, также имеет такой вид.
Это следует из
$$(a+br+cr^2+dr^3)(a-br+cr^2-dr^3)=(a+cr^2)^2-r^2(b+dr^2)^2\text{ и }(A+Br^2)(A-Br^2)=A^2-B^2r^4.$$

 \subsection*{Несколько извлечений квадратных корней }

  {\bf \ref{twoexa}.} Нет.

(a) Домножьте на сопряженное.

(b) Проще доказать сразу, что  $\sqrt[3]2\ne a+p\sqrt b+q\sqrt c+r\sqrt{bc}$, где $a,b,c,p,q,r\in\Q$.
Для этого достаточно доказать, что $\sqrt[3]2\ne u+v\sqrt c$, где
$u$ и $v$ --- числа вида $\alpha+\beta\sqrt b$ где $\alpha,\beta\in\Q$.
Идея доказательства в том, что числа такого вида $\alpha+\beta\sqrt b$ (с фиксированным $b$)
`ничуть не хуже' рациональных чисел. Т.е. сумма, разность, произведение
и частное чисел такого вида --- тоже число такого вида.
(Или, говоря научно, такие числа образуют {\it числовое поле}.)
Поэтому можно доказывать аналогично задаче \ref{number2}.e.
Ср. с задачей \ref{sopr4det}.a.

\smallskip
{\bf \ref{noncon}.}
(a) Предположим, что $\sqrt[3]2$ вещественно построимо.
Тогда существует такая башня квадратичных расширений
$$\Q=F_1\subset F_2\subset F_3\subset \ldots\subset F_{s-1}\subset F_s\subset\R,
\quad\mbox{что}\quad \sqrt[3]2\in F_r\setminus F_{r-1}.$$
Поскольку $\sqrt[3]2\not\in \Q$, то $r\ge3$.
Значит,
$$\sqrt[3]2=\alpha+\beta\sqrt a,\quad\mbox{где}\quad \alpha,\beta,a\in F_{r-1},
\quad \sqrt a\not\in F_{r-1}\quad \mbox{и} \quad \beta\ne0.$$
Отсюда
$$2=(\sqrt[3]2)^3=(\alpha^3+3\alpha\beta^2a)+(3\alpha^2\beta+\beta^3a)\sqrt a=u+v\sqrt a.$$
Поскольку $2\in\Q\subset F_{r-1}$, то $2-u\in F_{r-1}$.
Так как
$$v\sqrt a=2-u \quad \mbox{и} \quad v\in F_{r-1}, \quad \mbox{то} \quad 0=v=3\alpha^2\beta+\beta^3a.$$
Так как $3\alpha^2+\beta^2a>0$, получаем $\beta=0$ --- противоречие!

(b,c) Следует из \ref{tower2}.c.

(d) Обозначим $\varepsilon:=\varepsilon_7$.
Так как $\varepsilon\ne1$, то число $\varepsilon$ удовлетворяет уравнению 6-ой степени
$\varepsilon^6+\varepsilon^5+\varepsilon^4+\varepsilon^3+\varepsilon^2+ \varepsilon+1=0$.
Разделим обе части уравнения на $\varepsilon^3$.
Положим
$$f:=\varepsilon+\varepsilon^{-1},\quad\mbox{тогда}\quad \varepsilon^2+\varepsilon^{-2}=f^2-2\quad\mbox{и}
\quad \varepsilon^3+\varepsilon^{-3}=f(\varepsilon^2+\varepsilon^{-2}-1).$$
$$\text{Получим}\quad
f(f^2-3)+(f^2-2)+f+1=0,\quad\mbox{то есть}\quad f^3+f^2-2f-1=0.$$
Кандидаты на рациональные корни этого уравнения $f=\pm1$ отвергаются проверкой.
Значит, по теореме о кубических уравнениях (\ref{tower2}.c) число $f=\varepsilon+\varepsilon^{-1}$ не построимо.
Поэтому и $\varepsilon$ не построимо (поясните).

\smallskip
{\bf \ref{tower2}.} (a) Следует из (b).

(b) Это утверждение легко доказывается индукцией по количеству операций
калькулятора, необходимых для получения числа, с применением домножения на сопряженное.

(c) Часть `тогда' очевидна.
Чтобы доказать часть `только тогда', предположим, что хотя бы один из корней построим.
Для каждого из построимых корней $z$ рассмотрим минимальную цепочку расширений
$$\Q=Q_1\subset Q_2\subset Q_3\subset\ldots\subset Q_{r-1}\subset Q_r,\quad
\mbox{для которой}\quad z_1\in Q_r\setminus Q_{r-1}.$$
Возьмем корень $z=z_1$ с наименьшей длиной $l$ минимальной цепочки.

Если кубическое уравнение не имеет рациональных корней, то $l\ge2$.
Значит,
$$z_1=\alpha+\beta\sqrt a,\quad\mbox{где}\quad\alpha,\beta\in Q_{l-1},
\quad\sqrt a\not\in Q_{l-1}\quad \quad\mbox{и}\quad \beta\ne0.$$
Тогда по аналогу леммы о сопряжении
$z_2:=\overline z_1=\alpha-\beta\sqrt a$ также является корнем уравнения.
Поскольку
$$\beta\ne0,\quad\mbox{то}\quad \alpha-\beta\sqrt a\ne\alpha+\beta\sqrt a,
\quad\mbox{т. е.}\quad z_2\ne z_1.$$
Обозначим $z_3$ третий корень уравнения (возможно, $z_3\in\{z_1,z_2\}$).
По формуле Виета
$$z_1+z_2+z_3=(\alpha+\beta\sqrt a)+(\alpha-\beta\sqrt a)+z_3=2\alpha+z_3\in\Q,
\quad\mbox{поэтому}\quad z_3\in Q_{l-1}.$$
Следовательно, для корня $z_3$ существует цепочка меньшей длины, чем для $z_1$.
Противоречие.
QED


\smallskip
{\bf \ref{alge}.}
{\it Указание к (b).} Пусть $a=a_1$ и $b=b_1$ --- построимые числа, а $P$ и $Q$ ---
многочлены с рациональными коэффициентами минимальной степени, корнями которых
являются соответственно $a$ и $b$.
Пусть $a_2,\dots,a_m$ --- все остальные комплексные корни многочлена $P$, а $b_2,\dots,b_n$ --- все остальные комплексные корни многочлена $Q$.
Тогда
$$a+b,\quad a-b,\quad ab,\quad a/b,\quad \sqrt a \quad \text{--- корни многочленов}$$
$$P(x-b_1) \dots P(x-b_n),\quad P(x+b_1) \dots P(x+b_n),\quad P(x/b_1) \dots P(x/b_n),\quad P(xb_1) \dots P(xb_n),
\quad P(x^2),$$
соответственно.
Осталось доказать следующее утверждение.

Пусть $R(x,y)$ --- многочлен от двух переменных и $b_1,b_2,\dots,b_n$ ---
все комплексные корни другого многочлена.
Тогда многочлен $R(x,b_1) R(x,b_2) \dots R(x,b_n)$ от одной переменной
также имеет рациональные коэффициенты.

 \subsection*{К доказательству непостроимости в теореме Гаусса}

 {\bf \ref{irre}.}
Ответы (для старших коэффициентов 1):

(a) $P(x):=((x-\sqrt3)^2-2)((x+\sqrt3)^2-2)=(x^2-5)^2-24$; \quad

(b) $P(x-\sqrt5)P(x+\sqrt5)$; \quad (d) $(x^2-1)^2-3$; \quad (e) $((x^2-5)^2-1)^2-3$. \quad

(a) Для доказательства неприводимости примените лемму о сопряжении и получите, что каждое из 4 чисел $\pm\sqrt2\pm\sqrt3$ является корнем нужного многочлена.

(b) Аналогично (a).

(c) $\sqrt{2+\sqrt3}=\frac{1+\sqrt3}{\sqrt2}$.

(d) Аналогично (a) неприводимый многочлен с коэффициентами в $\Q[\sqrt 3]$ и корнем $\sqrt{1+\sqrt3}$
(и старшим коэффициентом 1) равен $x^2-2-\sqrt3$.
Значит, любой неприводимый многочлен $P$ с коэффициентами в $\Q$ и корнем $\sqrt{1+\sqrt3}$ делится на $x^2-2-\sqrt3$.
Применяя сопряжение относительно $\Q\subset\Q[\sqrt3]$ получаем, что $P$ делится и на $x^2-2+\sqrt3$.
Так как многочлены $x^2-2-\sqrt3$ и $x^2-2+\sqrt3$ взаимно просты, то $P$ делится на их произведение.

(e) Аналогично (d).

\smallskip
{\bf \ref{baslem}.} Аналогично \ref{irre}.de.
См. подробности в доказательстве непостроимости в теореме Гаусса в \S\ref{progau}.

\smallskip
{\bf \ref{irree}.}
Ответы (для старших коэффициентов 1): \quad
(5) $x^4+x^3+x^2+x+1$; \quad

(7) $x^6+x^5+\dots+x+1$; \quad (9) $x^6+x^3+1$; \quad (11) $x^{10}+x^9+\dots+x+1$; \quad
(25) $x^{20}+x^{15}+x^{10}+x^5+1$.

(5) Для доказательства неприводимости примените признак Эйзенштейна к многочлену $\Phi(x+1)=((x+1)^5-1)/x$ и
лемму Гаусса.

 \subsection*{Одно извлечение корня третьей степени}

  {\bf \ref{3ex}.} Нельзя.

(a) {\it Первое решение.}
Пусть можно. Тогда
$$3=(a^2+4bc)+(2ab+2c^2)\sqrt[3]2+(2ac+b^2)\sqrt[3]4.$$
Так как многочлен $x^3-2$ не имеет рациональных корней, то он неприводим над $\Q$.
Значит, $2ab+2c^2=2ac+b^2=0$ (ср. \ref{sopr3}.b).
Поэтому $b^3=-2abc=2c^3$.
Тогда либо $b=c=0$, либо $\sqrt[3]2=b/c$.
Оба случая невозможны.

Для остальных решений задачи \ref{3ex} обозначим $r=\sqrt[3]2$ и обозначим через $x_1,x_2,x_3\in\C$ числа,
заданные формулами (*), где $a,b,c\in\Q$.

{\it Второе решение.}
Пусть можно.
По лемме о сопряжении (\ref{sopr3}.a) многочлен $x^2-3$ имеет три корня $x_1,x_2,x_3$.
Так как ни один из них не рационален, то $b=c=0$ невозможно.
Значит, по лемме о линейной независимости (\ref{sopr}.c) эти корни различны.
Противоречие.

(b) Пусть можно.
Число $\cos(2\pi/9)$ является корнем уравнения $4x^3-3x=-\frac12$.
Два других его вещественных корня есть $\cos(8\pi/9)$ и $\cos(4\pi/9)$.
По лемме о сопряжении (\ref{sopr3}.a) многочлен $8x^3-6x-1$ имеет три корня $x_1,x_2,x_3$.
Так как ни один из них не рационален, то $b=c=0$ невозможно.
Значит, по лемме о линейной независимости (\ref{sopr3}.c) эти корни различны.

Заметим, что $\overline{\varepsilon^k}=\varepsilon^{-k}$.
Поэтому числа $a+br\varepsilon+cr^2\varepsilon^2$ и $a+br\varepsilon^2+cr^2\varepsilon$ комплексно сопряжены.
Значит, они не могут быть вещественными и различными.

(c) Пусть можно.
По лемме о сопряжении (\ref{sopr3}.a) многочлен $x^5-3$ имеет три корня $x_1,x_2,x_3$.
Так как ни один из корней не рационален, то $b=c=0$ невозможно.
Значит, по лемме о линейной независимости (\ref{sopr3}.c) эти корни различны.
Поэтому многочлен $x^5-3$ делится на $(x-x_1)(x-x_2)(x-x_3)$.
По лемме о рациональности (\ref{sopr3}.d) многочлен
\linebreak
$(x-x_1)(x-x_2)(x-x_3)$ имеет рациональные коэффициенты.
Противоречие с неприводимостью над $\Q$ многочлена $x^5-3$.

(d)  Аналогично (a), (b) получаем, что комплексные корни многочлена $x^3-3$ есть числа $x_1,x_2,x_3$.
Поэтому $(a+br+cr^2)\varepsilon^s=a+br\varepsilon+a+br\varepsilon^2$ для некоторого $s\in\{1,2\}$.
Отсюда по лемме о линейной независимости (\ref{sopr3}.c) $a=0$ и $bc=0$.
Поэтому либо $\sqrt[3]3=br$, либо $\sqrt[3]3=cr^2$.
Противоречие.

(e) Аналогично (b).

\smallskip
{\bf \ref{sopr3}.} (a) Аналогично \ref{sopr}.a и \ref{sopr4}.
Подставим $a+bt+ct^2$ в многочлен и поделим с остатком на $t^3-2$.
Подставляя $t=\sqrt[3]2$, получаем по лемме о линейной независимости (\ref{sopr3}.b), что остаток нулевой.
Значит, если $r^3=2$, то  $a+br+cr^2$ есть корень исходного многочлена.

(b) Так как многочлен $x^3-r^3$ не имеет рациональных корней, то он неприводим над $\Q$.

(c) Докажите отдельно для вещественной и мнимой части.

(d) {\it Первое решение} получается из тождества \ref{ident}.a подстановкой $x-a,br,cr^2$ вместо $a,b,c$, соответственно.

{\it Второе решение.} И при замене $r$ на $r\varepsilon$, и при замене $\varepsilon$ на $\varepsilon^2$ наш многочлен переходит в себя, даже если его рассматривать как многочлен от $x,r,\varepsilon$ по модулю $\varepsilon^3-1$.
Значит, его `коэффициент' при $x^kr^l$ при замене $\varepsilon$ на $\varepsilon^2$ переходит в себя.
Поэтому `коэффициент' при $x^kr^l$ не зависит от $\varepsilon$ и, следовательно, рационален.
Далее, его `коэффициент' при $x^k$ при замене $r$ на $r\varepsilon$  переходит в себя.
 Поэтому каждый `коэффициент' при $x^k$ рационален.
(Ввиду рассмотрения многочленов, а не чисел, не нужна даже лемма о линейной независимости.)


{\it Указание к третьему решению.}
Собирая вместе коэффициенты при одинаковых степенях $r$, имеем
$$x_1+x_2+x_3=3a+br(1+\varepsilon+\varepsilon^2)+cr^2(1+\varepsilon^2+\varepsilon)=3a\in\Q,$$
$$x_1^2+x_2^2+x_3^2=3a^2+2ab(1+\varepsilon+\varepsilon^2)r+(2ac+b^2)(1+\varepsilon^2+\varepsilon)r^2+6bcr^3+
c^2(1+\varepsilon+\varepsilon^2)r^4=3a^2+6bcr^3\in\Q$$
и $x_1^3+x_2^3+x_3^3=\dots$.

{\it Третье решение.}
Обозначим эти числа через $x_1,x_2,x_3$.
Ввиду теоремы Виета и формул Ньютона достаточно доказать рациональность чисел $x_1^k+x_2^k+x_3^k$ для $k=1,2,3$.
Обозначим $P(t):=a+bt+ct^2$.
Тогда $P(t)^k=a_0+a_1t+\dots+a_{2k}t^{2k}$ для некоторых $a_0,a_1,\dots,a_{2k}\in\Q$ (зависящих только от $k,a,b,c$).
Тогда
$$x_1^k+x_2^k+x_3^k=P(r)^k+P(\varepsilon r)^k+P(\varepsilon^2 r)^k = \sum\limits_{s=0}^{2k}a_sr^s(1+\varepsilon^s+\varepsilon^{2s}) =
3\sum\limits_{u=0}^{[2k/3]}a_{3u}r^{3u}\in\Q.$$
\quad

\smallskip
{\bf \ref{utv3}.} Аналогично задачам \ref{3ex}.abc.

\smallskip
{\bf \ref{que3}.} (3R,3C) Ответ:
многочлен $x^3+px+q$ имеет корень указанного вида тогда и только тогда, когда либо многочлен имеет рациональный
корень, либо

(3R) $(p/3)^3+(q/2)^2$ \qquad (3C) $|(p/3)^3+(q/2)^2|$

есть квадрат рационального числа. См. [Ak].

\smallskip
{\bf \ref{mno3}.} Пусть при извлечении корня третьей степени получилось число $r$.
Если $|r|\in\Q$, то утверждение очевидно.
Если $|r|\not\in\Q$, то многочлен $x^3-r^3$ неприводим над $\Q$.

Достаточно доказать, что $\frac1{a+br+cr^2}=h(r)$ для некоторого многочлена $h$.
Так как многочлены $x^3-r^3$ и $a+br+cr^2$ взаимно просты, то существуют многочлены $g$ и $h$,
для которых $h(x)(a+bx+cx^2)+g(x)(x^3-r^3)=1$.
Тогда $h$ --- искомый.

  \subsection*{Одно извлечение корня простой степени}

  {\bf \ref{kvadiffreal}.} Нельзя.

{\it Указание.} Используйте сформулированные ниже леммы.
Аналогично второму решению задачи \ref{3ex}.

Приведем решения. Обозначим $r:=\sqrt[7]2$ и $A(x):=a_0+a_1x+a_2x^2+\dots+a_6x^6$.

(a) Пусть можно.
Тогда по лемме о сопряжении (\ref{idea}.c) многочлен $x^2-3$ имеет корни $A(r\varepsilon_7^k)$ для
$k=0,1,2,\dots,6$.
По лемме о линейной независимости (\ref{idea}.e) они попарно различны.
Противоречие.

(b) Пусть можно.
Обозначим через $p$ многочлен, для которого $\cos7x=p(\cos x)$.
Тогда по леммам о сопряжении и о линейной независимости (\ref{idea}.c,e)
многочлен $2p(x)+1$ имеет попарно различные корни $A(r\varepsilon_7^k)$ для $k=0,1,2,\dots,6$.
Но все корни этого многочлена вещественны.


Имеем $\overline{\varepsilon_q^k}=\varepsilon_q^{-k}$.
Поэтому для $k>0$ числа $x_k:=A(r\varepsilon_q^k)$ и $x_{q-k}$ симметричны относительно вещественной оси (т.е. комплексно сопряжены).
Так как $q$ нечетно и числа $x_1,\dots,x_{q-1}$ попарно различны, то ни одно из них не может быть вещественным.

(c) Пусть можно.
Тогда по леммам о сопряжении и о линейной независимости (\ref{idea}.c,e) многочлен $x^{11}-3$ имеет попарно
различные корни $A(r\varepsilon_7^k)$ для $k=0,1,2,\dots,6$.
По лемме о рациональности (\ref{idea}.f) получаем противоречие с неприводимостью многочлена $x^{11}-3$.

(d) Аналогично (a), (b) получаем, что комплексные корни многочлена $x^7-3$ есть $A(r\varepsilon_7^k)$ для $k=0,1,2,\dots,6$.
Поэтому $A(r)\varepsilon_7^s=A(r\varepsilon_7)$ для некоторого $s\in\{1,2,3,4,5,6\}$.
Отсюда по лемме о линейной независимости (\ref{idea}.e) $a_k=0$ для любого $k\ne s$.
Поэтому $\sqrt[7]3=a_sr^s$.
Противоречие.

(e) Аналогично (b).

\smallskip
{\bf \ref{idea}.} (a) Все корни многочлена $x^q-r^q$ есть $r,r\varepsilon_q,r\varepsilon_q^2,\dots,r\varepsilon_q^{q-1}$.
Пусть он приводим над $\Q$.
Модуль свободного члена одного из сомножителей разложения раци
онален и равен произведению модулей некоторых $k$ из
этих корней, $0<k<q$.
Значит, $r^k\in\Q$.
Так как $q$ простое, то $kx+qy=1$ для некоторых целых $x,y$.
Тогда $r^{kx}=r(r^q)^{-y}$, откуда $r\in\Q$. Противоречие.

(b) Вытекает из (a).

(c) Аналогично задачам \ref{sopr}.a, \ref{sopr4} и \ref{sopr3}.a.
Используйте слабую лемму о линейной независимости (a).


(d) Пусть приводим. Аналогично доказательству неприводимости над $\Q$ получим $r\in\Q[\varepsilon_q]$.
Поэтому $r^2,r^3,\dots,r^{q-1}\in\Q[\varepsilon_q]$.
Составим таблицу $a_{kl}$ из рациональных чисел размера $q\times (q-1)$ из разложений чисел $r^k$ по степеням числа $\varepsilon_q$:
$$r^k=\sum\limits_{l=0}^{q-2}a_{kl}\varepsilon_q^l,\qquad 0\le k\le q-1.$$
При помощи прибавления к одной строке другой, умноженной на рациональное число, можно получить таблицу с нулевой строкой.
Значит, имеется многочлен степени меньше $q$ с корнем $r$.
Противоречие с неприводимостью многочлена $x^q-r^q$ над $\Q$.




(e) Вытекает из (d).

(f) Аналогично \ref{sopr3}.d, используя то, что при каждой из замен $r$ на $r\varepsilon_q$ и $\varepsilon_q$ на $\varepsilon_q^s$, $s=2,3,\dots,q-1$, наш многочлен переходит в себя.
Подробности приведены в доказательстве леммы о разложении в \S\ref{prorea}.

\smallskip
{\bf \ref{kronint1}.}
{\it Указание к (a).}
Предположим противное.
Аналогично  \ref{kvadiffreal}.a,b,c (и утверждению \ref{utv3}.a,b) получаем противоречие, используя леммы о сопряженности, о линейной независимости
и о рациональности (\ref{idea}.c,e,f).

{\it Решение (a).}
Предположим противное.
Тогда по леммам о сопряжении и о линейной независимости (\ref{idea}.c,e) данный многочлен $f$ имеет попарно различные корни $A(r\varepsilon_q^k)$ для $k=0,1,2,\dots,q-1$.
При $q>\deg f$ получаем противоречие.
При $q=\deg f$ получаем противоречие аналогично \ref{kvadiffreal}.b).
При $q<\deg f$ получаем противоречие по лемме о рациональности (\ref{idea}.f).

(b) Аналогично (a).
Для комплексного случая в леммах нужно использовать, что $r^q\ne b^q$ ни для какого $b\in\Q$.
В лемме о вещественности нужно начать со следующего рассуждения. Обозначим $x_k:=A(r\varepsilon_q^k)$.
Если заменить $r$ на $r\varepsilon_q$, то множество чисел $x_0,x_1,\dots,x_{q-1}$ не изменится, они лишь перенумеруются.
Поэтому можно считать, что $r\in\R$.

\smallskip
{\bf \ref{mnoq}.} Аналогично задачам \ref{mno3} и \ref{mno4}.

 \subsection*{Несколько извлечений корней}

 {\bf \ref{twocalc}.} (a) Нельзя.

(b)
Годится многочлен $f$, неприводимый и имеющий более одного вещественного корня.
Для последнего достаточно $f(0)>0$ и $f(1)<0$.
Например, можно взять $x^5-4x-2$.



\section{Доказательства неразрешимости  в радикалах}\label{pro}

\subsection{Лемма о калькуляторе и понятие поля}\label{profie}

Если $F\subset\C$, $r\in\C$ и $r^q\in F$ для некоторого целого положительного $q$, то обозначим
$$F[r]:=\{a_0+a_1r+a_2r^2+\dots+a_{q-1}r^{q-1}\ |\ a_0,\dots,a_{q-1}\in F\}.$$

{\bf Лемма о калькуляторе.} {\it Пусть $F\in\{\R,\C\}$.
Число $x\in F$ можно получить на $F$-калькуляторе тогда и только тогда, когда существуют такие
$r_1,\ldots r_{s-1}\in F$ и $q_1,\ldots q_{s-1}\in\Z$, что $q_k\ge 2$,
$$\Q=F_1\subset F_2\subset F_3\subset \ldots\subset F_{s-1}\subset F_s\ni \alpha,
\quad\mbox{где}\quad r_k^{q_k}\in F_k,\quad r_k\not\in F_k\quad\mbox{и}\quad F_{k+1}=F_k[r_k].$$
для любого $k=1,\dots,s-1$.}

\smallskip
Такая последовательность называется {\it башней расширений}.

{\it Полем} называется подмножество множества $\C$, замкнутое относительно операций сложения, умножения,
вычитания и деления на ненулевое число.
Это понятие полезно для нас тем, что теорема о делении с остатком и ее следствия верны для многочленов
с коэффициентами в поле.

Если $F$ поле, $q$ простое и $r\not\in F$, то многочлен $t^q-r^q$ неприводим над $F$
(аналогично слабой лемме о линейной независимости \ref{idea}.a, ср. с задачей \ref{mailemtow}).
Тогда $F[r]$ поле.

Напомним, что
$$\varepsilon_n:=\cos\dfrac{2\pi}n+i\sin\dfrac{2\pi}n.$$

\subsection{Доказательство непостроимости в теореме Гаусса}\label{progau}

Так как $\varepsilon_n=\varepsilon_{nk}^k$, то из построимости числа $\varepsilon_{nk}$ вытекает
построимость  числа $\varepsilon_n$.
Поэтому и по лемме о комплексификации для доказательства вещественной непостроимости в теореме Гаусса
достаточно показать, что  $\varepsilon_n$ непостроимо для

(A) $n$ простого, не представимого в виде $2^m+1$,

(B) $n=p^2$ квадрата простого.

\smallskip
{\bf Лемма о степенях двойки.} {\it Если неприводимый над $\Q$ многочлен $P$ с рациональными коэффициентами
имеет построимый корень, то $\deg P$ есть степень двойки.}

\smallskip
Эта лемма доказана далее.

Непостроимость числа $\varepsilon_n$ следует из леммы о калькуляторе и леммы о степенях двойки для
корня $\varepsilon_n$ многочлена

$\bullet$ $P(x):=x^{n-1}+x^{n-2}+\dots+x+1$ в случае (A) и

$\bullet$ $P(x):=x^{p(p-1)}+x^{p(p-2)}+\dots+x^p+1$ в случае (B).

Неприводимость этих многочленов $P(x)$ над $\Z$ вытекает из неприводимости многочленов $P(x+1)$ над $\Z$.
Последняя неприводимость доказывается применением следующего {\it признака Эйзенштейна}:

{\it Пусть $p$ простое. Если для многочлена с целыми коэффициентами старший коэффициент не делится на $p$,
остальные делятся на $p$, а свободный член не делится на $p^2$, то этот многочлен неприводим над $\Z$.}

Условие этого признака легко проверяется с помощью сравнения $(a+b)^p\equiv a^p+b^p\mod p$.

Неприводимость над $\Q$ вытекает из неприводимости над $\Z$ и следующей леммы Гаусса:

{\it Если многочлен с целыми коэффициентами неприводим над $\Z$, то он неприводим и над $\Q$.}

И признак Эйзенштейна, и  лемма Гаусса легко доказываются переходом к многочленам с коэффициентами $\Z_p$
(для леммы Гаусса рассмотрим разложение $P=P_1P_2$ данного полинома $P$ над $\Q$, возьмем целые $n_1$ и $n_2$
такие, что и $n_1P_1$, и $n_2P_2$, имеют целые коэффициеты, и возьмем простой делитель $p$ числа $n_1n_2$).

\smallskip
{\bf Лемма о сопряжении.} {\it Пусть $F\subset\C$ --- поле, $r\in\C-F$ и $r^2\in F$.
Определим отображение сопряжения
$\overline\cdot:F[r]\to F[r]\quad\mbox{формулой}\quad\overline{x+yr}:=x-yr.$
Это отображение корректно определено,}
$$\overline{z+w}=\overline z+\overline w, \quad \overline{zw}=\overline z\cdot\overline w\quad\mbox{and} \quad
\overline z=z\Leftrightarrow z=x+0r\in F.$$

Лемма о степенях двойки является случаем $k=1$ следующего утверждения.

\smallskip
{\bf Обобщенная лемма о степенях двойки.} {\it Если
$$\Q=F_1\subset F_2\subset F_3\subset \ldots\subset F_{s-1}\subset F_s\ni \alpha,
\quad\mbox{где}\quad r_k^2\in F_k,\quad r_k\not\in F_k\quad\mbox{и}\quad F_{k+1}=F_k[r_k]$$
для любого $k=1,\dots,s-1$,
то для каждого $k=1,2,\dots,s$ степень любого неприводимого над $F_k$ многочлена  с коэффициентами из $F_k$ и
корнем $\alpha$ есть степень двойки.}

\smallskip
{\it Доказательство.}
Индукция по $k$ вниз. База $k=s$ очевидна.
Докажем шаг.
Обозначим через $P_k$ любой неприводимый над $F_k$ многочлен с коэффициентами из $F_k$ и корнем $\alpha$.
Будем рассматривать делимость и НОД в $F_{k+1}$.
Так как $P_k(\alpha)=0$, то $P_k$ делится на $P_{k+1}$.
По лемме о сопряжении для $F=F_k$ и $F[r]=F_{k+1}$ имеем $P_k=\overline{P_k}$ делится на $\overline{P_{k+1}}$.
Обозначим $D:=GCD(P_{k+1},\overline{P_{k+1}})$.
Так как $P_{k+1}$ неприводим над $F_{k+1}$ и делится на $D$, то либо $D=1$, либо $P_{k+1}=D$.

Во втором случае так как $\overline D=D$, то $P_{k+1}=D\in F_k[x]$.
Значит, $P_k=P_{k+1}$ и шаг индукции доказан.

В первом случае $P_k$ делится на $M:=P_{k+1}\overline{P_{k+1}}$.
Так как $\overline M =M$, то $M\in F_k[x]$. Так как $P_k$ неприводим над $F_k$, то $P_k=M$.
Значит, $\deg P_k=2\deg P_{k+1}$ есть степень двойки по предположению индукции.
QED

\subsection{Доказательство неразрешимости в вещественных радикалах}\label{prorea}

{\bf Основная Лемма (вещественный случай).}
{\it Пусть $q$ простое, $F\subset\R$ поле, $r\in\R-F$ и $r^q\in F$.

(a) (линейная независимость)
Если $P(r)=0$ для некоторого многочлена $P\in F[\varepsilon_q][t]$ степени меньше $q$, то $P=0$.

(b) (сопряжение) Если $P\in F[\varepsilon_q][t]$ и $P(r)=0$, то $P(r\varepsilon_q^k)=0$ для любого $k=0,1,\dots,q-1$.}

\smallskip
{\it Доказательство части (a).}
Оба многочлена $P$ и $t^q-r^q$ с коэффициентами из $F[\varepsilon_q]$ имеют корень $r$.
Значит, их НОД имеет корень $r$ и степень $k$, $0<k\le\deg P<q$.
Все корни многочлена $t^q-r^q$ есть $r,r\varepsilon_q,r\varepsilon_q^2,\dots,r\varepsilon_q^{q-1}$.
Свободный член НОД'а равен произведению некоторых $k$ из этих корней.
Тогда $r^k\in F[\varepsilon_q]$.
Так как $q$ простое, то $kx+qy=1$ для некоторых целых $x,y$.
Тогда $r = (r^k)^x (r^q)^y\in F[\varepsilon_q]$.
\footnote{Другая запись следующего абзаца с использованием понятия размерности:
тогда $\dim_F F[r]\le \dim_F F[\varepsilon_q]\le q-1$.}

Поэтому $r^2,r^3,\dots,r^{q-1}\in F[\varepsilon_q]$.
Составим таблицу $a_{kl}\in F$ размера $q\times (q-1)$ из разложений чисел $r^k$ по степеням числа $\varepsilon_q$:
$$r^k=\sum\limits_{l=0}^{q-2}a_{kl}\varepsilon_q^l,\qquad 0\le k\le q-1.$$
При помощи нескольких операций прибавления к одной строке другой, умноженной на число из $F$,
можно получить таблицу с нулевой строкой.

Значит, имеется ненулевой многочлен $P_1$ степени меньше $q$ с коэффициентами из $F$ и корнем $r$.
Дальнейшие рассуждения аналогичны первому абзацу.
Нужно только заменить $P$ на $P_1$ и $F[\varepsilon_q]$ на $F$.
Получаем $r\in F$ --- противоречие.
QED

\smallskip
{\it Доказательство части (b).}
Так как $P(r)=0$, то остаток от деления многочлена $P(t)$ на $t^q-r^q$ принимает значение 0 в точке $r$.
Значит, по (a) этот остаток равен нулю.
Отсюда вытекает заключение части (b).
QED

\smallskip
{\it Доказательство теоремы о неразрешимости в вещественных радикалах.}
Предположим, напротив, что некоторый корень $x_0$ уравнения $8x^3-6x+1=0$ можно получить на вещественном калькуляторе.
Тогда по лемме о калькуляторе для $F=\R$ существует наименьшее $s$, для которого найдется башня расширений,
последнее поле $F_s$ которой содержит некоторый корень $x_1$ уравнения $8x^3-6x+1=0$ (возможно, $x_1\ne x_0$).
Обозначим $F:=F_{s-1}$, $q:=q_{s-1}$ и $r:=r_{s-1}$.
Тогда $x_1=h(r)$ для некоторого многочлена $h$  с коэффициентами в $F$ степени больше 0 и меньше $q$.

Применим основную леммe (вещественный случай) (b) к многочлену $P(t):=8h(t)^3-6h(t)+1$.
Так как $8h(r)^3-6h(r)+1=0$, получим, что $h(r\varepsilon_q^k)$ является корнем уравнения $8x^3-6x+1=0$ для любого $k=0,1,\dots,q-1$.
Если $h(r\varepsilon_q^k)=h(r\varepsilon_q^l)$ для некоторых $0\le k<l\le q-1$, то по основную леммe (вещественный случай) (a) получим $\deg h=0$ --- противоречие с минимальностью $s$.

Итак, числа $h(r\varepsilon_q^k)$, $0\le k\le q-1$, --- попарно различные корни уравнения $8x^3-6x+1=0$.
Значит, $q=2$ или $q=3$.
Если $q=2$, то по теореме Виета третий корень уравнения $8x^3-6x+1=0$ равен $-2h(0)\in F$ --- противоречие
с минимальностью числа $s$.
Если $q=3$, то из $\overline{\varepsilon_3}=\varepsilon_3^2$ вытекает $\overline{h(r\varepsilon_3)}=h(r\varepsilon_3^2)$.
Это противоречит вещественности и различности последних двух чисел.
QED

\subsection{Доказательство неразрешимости в радикалах}\label{progal}

Теорема Галуа (о неразрешимости в радикалах) вытекает из следующего результата.
Он интересен и нетривиален даже для многочленов пятой степени.

\smallskip
{\bf Теорема Кронекера.} {\it Если многочлен простой степени неприводим над $\Q$, имеет более одного вещественного корня и хотя бы один невещественный, то ни один из его корней невозможно получить на комплексном калькуляторе.}

\smallskip
Из следующих лемм только основная лемма (a,c) и лемма об уплотнении прямо используются в доказательстве
теоремы Кронекера.
Основная лемма (b) используется только для основной леммы (c) .

Основная лемма (c) показывает, что при некоторых предположениях все корни многочлена $g\in F[x]$
являются значениями в `сопряженных' точках некоторого многочлена из $F[x]$.

\smallskip
{\bf Основная Лемма.}
{\it Пусть $q$ простое, $F\subset\C$ поле, $r\in\C-F$ и $r^q,\varepsilon_q\in F$.

(a) (линейная независимость) Если $P(r)=0$ для некоторого многочлена $P\in F[t]$ степени меньше $q$, то $P=0$.
\footnote{
Аналог леммы о линейной независимости без условия `$\varepsilon_q\in F$' неверен для $q>2$, $F=\R$ и $r=\varepsilon_q$.
Условие `$\varepsilon_q\in F$'
пропущено в [P, стр. 580-581].
В [P]
утверждение `$q=p$' вверху стр. 581 (для $p=2$)
означает следующее: если квадратный трехчлен $f$ неприводим над полем $k$, содержащим $i$,
и приводим над $k[\sqrt[q]a]$ для некоторого $a\in K$ и простого $q$, то $q=2$.
Это неверно для $f(x)=x^2+x+1$, $q=3$, $a=1$ и $k=\Q[i]$.
Ошибка в доказательстве в [P] --- в предыдущем предложении: (верную)
теорему 1 на стр. 572 применить нельзя, т.к. возможно $a=b^q$ для некоторого $b\in k$ (хоть $\sqrt[q]a\not\in k$).
}

Или, эквивалентно, для любого $\alpha\in F[\varepsilon_q,r]$ существуют единственные
$a_0,a_1,\dots,a_{q-1}\in F[\varepsilon_q]$, для которых $x=a_0+a_1r+a_2r^2+\dots+a_{q-1}r^{q-1}$.

(b) (сопряжение) Если $P\in F[x,t]$ и $P(x,r)=0$ как многочлен от $x$, то $P(x,r\varepsilon_q^k)=0$
как многочлен от $x$ для любого $k=0,1,\dots,q-1$.

(c) (разложение) Если многочлен $g$ простой степени с коэффициентами в $F$ неприводим над $F$ и приводим над $F[r]$,
то $g(x)=A(x-x_0)(x-x_1)\dots (x-x_{q-1})$ для некоторого $A\in F$ и попарно различных значений $x_0,x_1,\dots,x_{q-1}$ некоторого многочлена с коэффициентами из $F$ в точках $r,r\varepsilon_q,\dots,r\varepsilon_q^{q-1}$.}

\smallskip
{\it Доказательство части (a).}
Аналогично первому абзацу доказательства вещественной леммы о линейной независимости, с заменой
$F[\varepsilon_q]$ на $F$.
QED



\smallskip
{\it Доказательство части (b).}
Утверждение части (b) инвариантно относительно деления многочлена $P$ с остатком на $t^q-r^q$.
Поэтому можно считать, что $\deg_tP<q$.
В этом случае часть (b) получается покоэффициентным применением части (a).
QED

\smallskip
{\it Доказательство части (c).}
По условию существует приведенный неприводимый над $F[r]$ делитель многочлена $g$ в $F[r]$.
Этот делитель есть значение $h(x,r)$ в точке $r$ некоторого многочлена $h\in F[x,t]$ степени по $t$ больше 0,
а по $x$ меньше $\deg g$.
Итак, $h(x,r)$ неприводим над $F[r]$ и $g(x)=h(x,r)h_1(x,r)$ для некоторого многочлена $h_1\in F[x,t]$.
Обозначим $\varepsilon:=\varepsilon_q$.

Применим (b) к $P(x,r):=g(x)-h(x,r)h_1(x,r)$.
Получим, что $g(x)$ делится на многочлен  $h(x,r\varepsilon^k)$ в $F[r]$ для любого $k=0,1,\dots,q-1$.

Многочлен  $h(x,r\varepsilon^k)$ неприводим над $F[r]$ для любого $k=0,1,\dots,q-1$.

(Иначе применим
(b) к многочлену $P$, равному разности $h(x,r\varepsilon^k)$ и его сомножителей.
Получим, что многочлен $h(x,r)$ приводим над $F[r]$.
Противоречие.)

По (a) многочлены $h(x,r\varepsilon^k)$ различны для различных $k=0,1,\dots,q-1$.
Значит, $g$ делится на их произведение.
По (a) это произведение можно однозначно представить в виде
$$a_0(x)+a_1(x)r+\dots+a_{q-1}(x)r^{q-1}\quad\text{для некоторых}\quad a_k\in F[x].$$
Так как произведение переходит в себя при заменe $r\to r\varepsilon$ (которая корректно определена по (a)),
то по (a) $a_k(x)=a_k(x)\varepsilon^k\in F[x]$ для любого $k=1,2,\dots,q-1$.
Отсюда $a_k(x)=0$ для любого $k=1,2,\dots,q-1$.
Значит, произведение равно $a_0(x)\in F[x]$.

Из этого и неприводимости $g$ над $F$ следует, что $g$ равно этому произведению.
Тогда $\deg g=q\deg_xh$.
Так как $\deg g$ простое и $\deg_xh<\deg g$, то $\deg_xh=1$ (и $\deg g=q$).
Значит, $-h(0,r)\in F[r]$ есть корень многочлена $g$.
И остальные его корни есть $-h(0,r\varepsilon^k)$ для $k=0,1,\dots,q-1$.
QED

\smallskip
{\bf Лемма об уплотнении башни расширений.}
{\it Если число можно получить на комплексном калькуляторе, то существует башня расширений из леммы о калькуляторе для $F=\C$,
для которой при любом $k=1,2,\dots,s-1$ число $q_k$ простое, $\varepsilon_{q_k}\in F_k$ и
либо $r_k\in\R$, либо $|r_k|^2\in F_k$.}

\smallskip
{\it Доказательство.}
При помощи индукции `вниз' по $q$ покажем, что из произвольной башни расширений можно получить башню расширений,
для которой $\varepsilon_{q_k}\in F_k$ при любом $q_k>q$.
Тогда при $q=1$ получим башню расширений, для которой $\varepsilon_{q_k}\in F_k$ при любом $k=1,2,\dots,s-1$.
База: $q=\max_k q_k$; в этом случае доказывать нечего.
Для доказательства шага индукции возьмем наименьшее такое $k$, что $q_k=q$.
Вставим между $F_{k-1}$ и $F_k$ `получение $\varepsilon_q$ при помощи корней степени меньше $q$'
из теоремы Гаусса о понижении \ref{lower}.
Если такого $k$ нет, то шаг индукции очевиден.

Далее заменим извлечение корня составной степени $ab$ на пару извлечений корней $a$-й и $b$-й степеней.
Условие $\varepsilon_{q_k}\in F_k$ при любом $k=1,2,\dots,s-1$ сохранится, ибо если $\varepsilon_{ab}\in F_k$,
то $\varepsilon_a\in F_k$ и $\varepsilon_b\in F_k$.

Назовем башню расширений {\it интересной}, если при любом $k=1,2,\dots,s-1$ число $q_k$ простое и
$\varepsilon_{q_k}\in F_k$.
При помощи индукции `вниз' по $l$ покажем, что из произвольной интересной башни можно получить интересную башню,
для которой при любом $k\le s-l$
$$\overline F_k=F_k\quad\text{и \quad либо}\quad r_k\in\R,\quad\text{либо}\quad |r_k|^2\in F_k.$$
Тогда при $l=0$ получим утверждение леммы.
База: $l=s-1$; в этом случае доказывать нечего.
Докажем шаг индукции.
(Если $r_k\in\R$, то шаг индукции очевиден, но следующее рассуждение тоже проходит.)
Так как $F_k=\overline{F_k}$ и $r_k^{q_k}\in F_k$, то $|r_k|^{2q_k}=r_k^{q_k}\overline{r_k^{q_k}}\in F_k$.
Поэтому $F_k[|r_k|^2]=F_k[\sqrt[q_k]{|r_k|^{2q_k}}]$, где берется вещественное значение корня.
Заменим подбашню
$$F_k\subset F_{k+1}\subset\dots\subset F_s\quad\text{на подбашню}\quad F_k\subset F_k[|r_k|^2]
\subset F_k[r_k,\overline r_k]=F_{k+1}[\overline r_k]\subset\dots \subset F_s[\overline r_k].$$
Ясно, что интересность башни сохраняется при этой замене.
Далее в новой подбашне из каждого набора совпадающих соседних полей оставляем только одно.
После чего пользуемся предположением индукции.
QED

\smallskip
{\it Доказательство теоремы Кронекера.}
Предположим, напротив, что некоторый корень данного многочлена $g$ можно получить на комплексном калькуляторе.
Тогда возьмем башню расширений из леммы об уплотнении башни расширений.
Так как $g$ неприводим над $\Q$ и приводим над последним полем башни, то существует такое $s$, что $g$ неприводим над $F_s$ и приводим над $F_{s+1}$.
Обозначим $r:=r_s$ и $q:=q_s$.
По основной лемме (c) $g(x)=A(x-x_0)(x-x_1)\dots (x-x_{q-1})$ для некоторого $A\in F_s$ и различных значений $x_0,x_1,\dots,x_{q-1}$ некоторого многочлена $a_0+a_1t+\dots+a_{q-1}t^{q-1}$ с коэффициентами из $F_s$ в точках $r\varepsilon_q^k$, $0\le k\le q-1$.
Вещественность значения $x_k$ равносильна тому, что $x_k=\overline x_k$.
Заметим, что $\overline{\varepsilon_q^k}=\varepsilon_q^{-k}$.

Если $r\in\R$, то по основной лемме (a) для любого $k\in\{0,1,\dots,q-1\}$ условие $x_k=\overline x_k$ равносильно
тому, что $a_s\varepsilon^{2sk}=\overline a_s$ для любого $s=0,1,\dots,q-1$.
Следовательно, $x_k\in\R$ не более, чем для одного $k$.

Если $r\not\in\R$, то лемме об уплотнении $|r|^2\in F_s$.
 Тогда  $\overline r^s=\dfrac{|r|^{2s}}{r^q}r^{q-s}$, где $\dfrac{|r|^{2s}}{r^q}\in F_s$.
Значит, по основной лемме (a) для любого $k\in\{0,1,\dots,q-1\}$ условие $x_k=\overline x_k$ равносильно
тому, что $a_0=\overline a_0$ и $a_s=\overline a_{q-s}\dfrac{|r|^{2q-2s}}{r^q}$ для любого $s=1,2,\dots,q-1$.
Эти равенства не зависят от $k$.
Поэтому если среди чисел $x_0,\dots,x_{q-1}$ есть вещественное, то все они вещественны.

Противоречие.
QED

\smallskip
{\bf Замечания.}
Отличие доказательства теоремы Кронекера от доказательства теоремы неразрешимости в вещественных радикалах заключается

$\bullet$ в `комплексификации': $r^q$ и коэффициенты многочлена $h$ могут быть комплексными.

$\bullet$ в необходимости доказывать наличие корня у многочлена, неприводимого над $F_{s-1}$ и приводимого над $F_s$.
(Если брать наименьшее $s$ такое, что данный многочлен имеет корень в $F_s$, то нужно доказывать неприводимость над $F_{s-1}$;  это менее удобно.)

Наличие корня --- нетривиальная часть основной леммы (c) (разложение из наличия корня вывести легко).
Чтобы доказывать основную лемму (с), а не ее усиленную версию из следующего пункта (и, тем самым, обойтись без теоремы о размерности башни) в лемме об уплотнении мы добиваемся условия $\varepsilon_{q_k}\in F_k$.
Для сильной теоремы о неразрешимости в вещественных радикалах этот трюк не проходит, поэтому ее доказательство более сложно (в частности, использует теорему о размерности башни).


Два разбираемых случая в конце доказательства теоремы Кронекера немного отличаются от случаев, разобранных в конце доказательства из [T].

В начале 2-й колонки на стр. 14 в [T] фактически используется, что $\rho\in R$.
А это неверно без дополнительных стараний типа леммы об уплотнении башни расширений.


  \subsection{Сильная вещественная теорема о неразрешимости}\label{prostr}

 {\bf Сильная вещественная теорема о неразрешимости.}
{\it Если многочлен простой нечетной степени неприводим над $\Q$ и имеет более одного вещественного корня,
то ни один из его корней невозможно получить на вещественном калькуляторе.}

\smallskip
{\bf Сильная лемма о разложении.}
{\it Пусть $F\subset\C$ поле, $q$ простое, $r\in\C$, \ $r^q\in F$.
Если многочлен $g$ простой степени с коэффициентами в $F$ неприводим над $F$ и приводим над $F[r]$, то
$g(x)=A(x-x_0)(x-x_1)\dots (x-x_{q-1})$ для некоторого $A\in F$ и попарно различных значений $x_0,x_1,\dots,x_{q-1}$ некоторого многочлена с коэффициентами из $F[\varepsilon_q]$ в точках $r,r\varepsilon_q,\dots,r\varepsilon_q^{q-1}$.}

\smallskip
Лемма интересна даже для $F\subset\R$ и даже для $F=\Q$, хотя неразрешимость за одно извлечение корня доказывается и без нее.


\smallskip
{\bf Лемма о потере неприводимости.}
{\it Если многочлен $g$ простой степени неприводим над полем $F$ и приводим над $F[\varepsilon_q]$, то $q>\deg g$.}

\smallskip
Доказательство леммы приводится в конце этого пункта.

\smallskip
{\it Доказательство сильной леммы о разложении.}
Аналогично доказательству основной леммы (c).
Везде, кроме последнего абзаца, нужно заменить $F$ на $F[\varepsilon_q]$.
Перед последним абзацем нужно ставить следующее:
{\it `Так как $g$ делится на произведение, то $\deg g\ge q$.
Из этого и леммы о потере неприводимости следует, что $g$ неприводим над $F[\varepsilon_q]$. '}
QED

\smallskip
{\it Доказательство сильной вещественной теоремы о неразрешимости.}
Аналогично доказательству теоремы Кронекера.
Отличие только в том, что вместо леммы об уплотнении мы используем только простоту всех $q_k$, полагаем $F:=F_{s-1}[\varepsilon_q]$, вместо основной леммы (c) используем сильную лемму о разложении
и не рассматриваем случай $r\not\in\R$.
QED



\smallskip
Осталось доказать лемму о потере неприводимости.
\footnote{Было бы интересно доказать эту лемму (а значит, и сильную вещественную теорему о неразрешимости)
без использования теоремы о размерности башни.}

 Для полей $K\subset L$ {\it размерностью $\dim_KL$ поля $L$ над полем $K$} называется наименьшее $s$, для которого существуют $s$ таких элементов $l_1,\dots,l_s\in L$, что для любого $l\in L$ существуют $k_1,\dots,k_s\in K$, для которых $l=k_1l_1+\dots+k_sl_s$.

\begin{pr}\label{dim}
(a) Найдите $\dim_{\Q}\Q[\sqrt[5]3]$.

(b) Если $\beta$ --- корень неприводимого над полем $F$ многочлена $g$, то $\dim_FF[\beta]=\deg g$.

(c) {\bf Теорема о размерности башни.} Для любых полей $K\subset L\subset M$  выполнено $\dim_K M=\dim_LM\cdot\dim_KL$.
\end{pr}

{\it Доказательство леммы о потере неприводимости.}
 Обозначим через $\beta\in \C$ произвольный корень многочлена $g$.
По теореме о размерности башни
$$\dim_{F[\varepsilon_q]}F[\beta,\varepsilon_q]\cdot\dim_FF[\varepsilon_q]=\dim_{F[\beta]}F[\beta,\varepsilon_q]\cdot
\dim_FF[\beta].$$
Так как $g$ неприводим над $F$, то $\dim_FF[\beta]=\deg g$.
Так как $g$ приводим над $F[\varepsilon_q]$, то $\dim_{F[\varepsilon_q]}F[\beta,\varepsilon_q]<\deg g$.
Так как $\deg g$ простое, то $\dim_FF[\varepsilon_q]$ делится на $\deg g$.
Так как $\dim_FF[\varepsilon_q]<q$, то $q>\deg g$.
(На самом деле, даже $q-1>\deg g$.)
 QED

\bigskip
{\bf Литература}

[A] В.Б. Алексеев, Теорема Абеля. М: Наука, 1976.

[Ak] Д. Ахтямов, Решение кубических уравнений при помощи одного извлечения корня,
доклад на ММКШ-2013, http://www.mccme.ru/circles/oim/mmks/works2013/akhtyamov2.pdf

[B] J. Bergen, A Concrete Approach to Abstract Algebra: From the Integers to the
Insolvability of the Quintic, 2010.
http://www.elsevierdirect.com/product.jsp?isbn=9780123749413

[BK] Бурда Ю., Кадец Л. Семнадцатиугольник и закон взаимности Гаусса, Мат. Просвещение, 17 (2013).
http://www.mccme.ru/free-books/matprosi.html.

[Ch] Н. Н. Чеботарев, Основы теории Галуа. Часть 1. Л., М.: Гостехиздат, 1934.

[Ch1] Г.Р. Челноков, Основы теории Галуа в интересных задачах,
\linebreak
http://www.mccme.ru/circles/oim/materials/grishalois.pdf.
(засмотрено 11.11.2010)

[CR] Р. Курант, Дж. Роббинс, Что такое математика. М.: МЦНМО, 2004.

[D] H. D\"orrie, 100 Great Problems of Elementary Mathematics: Their History and Solution, Dover Publ, New York, 1965.

[E1] H.M. Edwards, Galois Theory, Springer Verlag, 1984.

[E2] H.M. Edwards, The construction of solvable polynomials,
Bull. Amer. Math. Soc. 46 (2009), 397-411.
Errata: Bull. Amer. Math. Soc. 46 (2009), 703-704.

http://www.ams.org/journals/bull/2009-46-03/S0273-0979-09-01253-1/

[FT] D. Fuchs, S. Tabachnikov,  Mathematical Omnibus. AMS, 2007.
\linebreak
http://www.math.psu.edu/tabachni/Books/taba.pdf

[Ga] К. Ф. Гаусс, Арифметические исследования. Труды по теории чисел.
М.: Изд-во АН СССР, 1959. С.~9--580.


[Gi] С. Гиндикин, Дебют Гаусса, Квант, 1972 N1, 2--11.

[H] Ch.R. Hadlock, Field Theory and its Classical Problems, Carus Mathematical Monographs 19,
The Mathematical Association of America, 1978.
\linebreak
http://books.google.com/books?id=5s1p0CyafnEC\&printsec=frontcover\&dq=hadlock.

[K] А.А. Кириллов, О правильных многоугольниках, функции Эйлера и числах Ферма,
Квант, 1977 N7, 2--9 или Квант, 1994 N6, 15--18.

[Ki] В.А. Кириченко, Построения циркулем и линейкой и теория Галуа,
\linebreak
http://www.mccme.ru//dubna/2005/courses/kirichenko.html

[Ko] В.А. Колосов, Теоремы и задачи алгебры, теории чисел и комбинаторики.
М: Гелиос, 2001.

[Kh1] А.Г. Хованский, Топологическая теория Галуа, Москва, МЦНМО, 200?

[Kh2] А. Г. Хованский Построения циркулем и линейкой, Мат. Просвещение, 17 (2013).
http://www.mccme.ru/free-books/matprosi.html.

[KS] П. Козлов и А. Скопенков, В поисках утраченной алгебры: в направлении
Гаусса (подборка задач), Мат. Просвещение, 12 (2008), 127--144,
http://arxiv.org/abs/0804.4357 (v1)

[Le] L. Lerner, Galois Theory without abstract algebra, http://arxiv.org/abs/1108.4593.

[Li] Дж. Литлвуд, Математическая смесь. М.: Наука,   1978.

[Ma] Ю. И. Манин, О разрешимости задач на построение с помощью циркуля и
линейки.
В кн. Энциклопедия элементарной математики. Книга четвертая (геометрия).
Под редакцией П. С. Александрова, А. И. Маркушевича и А. Я. Хинчина.
М., Физматгиз, 1963.

[Po] М. М. Постников, Теория Галуа. М.: Гос. изд-во физ.-мат. л-ры, 1963.

[P] В.В. Прасолов,
Задачи по алгебре, арифметике и анализу (М.: МЦНМО, 2007)
\linebreak
ftp://ftp.mccme.ru/users/prasolov/algebra/algebra2.pdf


[PS] В.В. Прасолов и Ю.П. Соловьев,
Эллиптические функции и алгебраические уравнения. М.: Факториал, 1997.
http://www.mccme.ru/prasolov

[S] A. Skopenkov, A simple proof of the Abel-Ruffini theorem,
Mat. Prosveschenie, 15 (2011) 113-126,
http://arxiv.org/abs/1102.2100.



[Sa] А. Сафин, Программа для построения правильных многоугольников циркулем и линейкой, http://www.mccme.ru/mmks/dec08/Safin.pdf

[Sa1] Теорема Кронекера, http://www.cdoosh.ru/lmsh/archive.html, 2011, 10 класс.

[T] В.М. Тихомиров, Абель и его великая теорема, Квант, 2003, N1.
\linebreak
http://kvant.mccme.ru/pdf/2003/01/kv0103abel.pdf

[Va] Вагутен Н., Сопряженные числа. Квант, 1980, N2,
\linebreak
http://kvant.mccme.ru/1980/02/sopryazhennye\_chisla.htm


[Vi] Э. Б. Винберг, Алгебра многочленов. М.: Просвещение, 1980.

[W] Б.Л. ван дер Варден, Алгебра, М: Наука, 1976.


\end{document}